\documentclass[12pt]{amsart}

\usepackage[margin=1in]{geometry}
\usepackage[metapost]{mfpic}
\usepackage{indentfirst,latexsym,amsmath,amssymb,amsthm}
\usepackage{enumerate}
\usepackage{mathptmx}
\usepackage{amsthm}
\usepackage{hyperref}
\usepackage{cleveref}
\usepackage{geometry}

\usepackage{t1enc}
\usepackage[cp1250]{inputenc}

\usepackage{mathrsfs}
\usepackage{calc}

\usepackage{color}
\usepackage{colortbl}
\definecolor{labelkey}{rgb}{0,0,1}

\newtheorem{thm}{Theorem}[section]
\newtheorem{prop}[thm]{Proposition}
\newtheorem{lm}[thm]{Lemma}

 \newtheorem{cor}[thm]{Corollary}
 \newtheorem{rem}[thm]{Remark}

\numberwithin{equation}{section}

\def\N{\mathbb{N}}

 \newcommand{\eps}{\varepsilon}
\newcommand{\ov}{\overline}

\renewcommand{\leq}{\leqslant}
\renewcommand{\geq}{\geqslant}

\def\vp{\varphi}
\def\A{\mathbf{A}}
\def\F{\mathbf{F}}
\def\X{\mathbf{X}}
\def\u{\mathbf{u}}
\def\w{\mathbf{w}}
\def\K{\mathbf{K}}
\def\lma{\lambda}

\def\Gr{\mathrm{Gr}}
\def\Epi{\mathrm{Epi\,}}
\def\R{\mathbb{R}}
\def\wew{\mathrm{int\,}}
\def\part{\partial}
\def\wt{\widetilde}
\opengraphsfile{myfigs}

\def\la{\langle}
\def\ra{\rangle}

\author[A. \'Cwiszewski]{Aleksander \'Cwiszewski}
\author[G. Gabor]{Grzegorz Gabor}
\address{Faculty of Mathematics and Computer Science. Nicolaus Copernicus university in Toru\'n, Poland}
\email{aleks@mat.umk.pl, ggabor@mat.umk.pl}
\author[W. Kryszewski]{Wojciech Kryszewski}
\address{Institute of Mathematics, Lodz University of Technology, Lodz, Poland}
\email{wojciech.kryszewski@p.lodz.pl}

\title[Invariance and Strict invariance for Nonlinear Evolutions Problems]{Invariance and Strict Invariance for Nonlinear Evolution Problems with Applications}

\date{\today}

\subjclass[2010]{37L05 (Primary) 47H06, 47J35, 35K91 (Secondary)}
\keywords{Evolution equation, invariance, strict invariance, PDEs, accretive operator}
\numberwithin{equation}{section}

\begin{document}
\pagestyle{myheadings}
\baselineskip15pt

\begin{abstract}
Sufficient conditions for the invariance of evolution problems governed by perturbations of (possibly nonlinear) $m$-accretive operators are provided. The conditions for the invariance with respect to sublevel sets of a constraint functional are expressed in terms of the Dini derivative of that functional, outside the considered sublevel set in directions determined by the governing $m$-accretive operator. An approach for non-reflexive Banach spaces is developed and some result improving a recent paper \cite{CaDaPraFran} is presented. Applications to nonlinear obstacle problems and age-structured population models are presented in spaces of continuous functions where advantages of that approach are taken. Moreover, some new abstract criteria for the so-called strict invariance are derived and their direct applications to problems with barriers are shown.
\end{abstract}

\maketitle
\vspace{-7mm}

\section{Introduction}
In the paper we study evolution problems of the form
\begin{equation}\label{equat-A+f}
\left\{
\begin{array}{ll}
\dot u \in -Au +f(u),\; t\geq 0,\\
u(0)=x_0\in \overline{D(A)}\cap \Omega,
\end{array}\right.
\end{equation}
where $A:D(A)\multimap X$ is a quasi-$m$-accretive operator in a Banach space $(X,\|\cdot\|)$ and $f:\Omega\to X$, with open $\Omega\subset X$, is a continuous map.
We are interested in the behavior of the so-called {\em integral} solutions $u$ to \cref{equat-A+f} related to the closed set $K\subset \Omega\cap \overline{D(A)}$.
In particular we look for conditions implying the {\em invariance} of $K$ with respect to the evolution determined by \eqref{equat-A+f}, meaning that {\em any} integral solution $u$ to \eqref{equat-A+f} starting at $x_0\in K$ remains in $K$, i.e. $u(t)\in K$ for every $t$ from the maximal interval of existence $[0,\tau_u)$, $0<\tau_u\leq\infty$.  Moreover issues concerning {\em strict invariance}, when {\em all} solutions starting at $x_0\in K$ stay for $t\in (0,\tau_u)$ in the interior $\wew K$ of $K$ will be of interest. Note that the concepts of invariance and strict invariance differ from the so-called  {\em viablity} of $K$ with respect to \eqref{equat-A+f} meaning that for any $x_0\in K$ {\em there exists} a solution $u$ starting at $x_0$ and staying in $K$ (see, e.g. \cite{au-viab}, \cite{CNV} and the extensive bibliography therein).\\
\indent  Let us start a discussion with some particular results and assume that $A\equiv 0$, i.e. consider the problem
\begin{equation}\label{equat+f}
\left\{
\begin{array}{ll}
\dot u = f(u),\; t\geq 0,\\
u(0)=x_0\in X,
\end{array}\right.
\end{equation}
(in this case integral and $C^1$-solutions to \eqref{equat+f} coincide) and let $K\subset \Omega$ be closed. In view of, e.g., \cite[Theorems 4.1.2,  4.1.3]{CNV} $K$ is invariant with respect to \eqref{equat+f} provided it is locally viable and $f$ satisfies a one-sided estimate of the following form
\begin{equation}\label{one sided}\big[x-y,f(x)-f(y)\big]_+\leq \omega\big(\|x-y\|\big)\; \text{for}\; x\in U\setminus K\; \text{and}\; y\in K,\end{equation}
where $U\subset\Omega$ is a neighborhood of $K$, $\omega$ is a {\em uniqueness function} (see Subsection \ref{uf}) and $[\cdot,\cdot]_+$ stands for the {\em right semi-inner product} (see Subsection \ref{semi}). In the proof condition \eqref{one sided} is used to compare an arbitrary solution $u$ of \eqref{equat+f} starting at some $x_0\in K$ and leaving $K$ with a viable local solution $v$ starting at $x_0\in K$ and then to  show  that $u$ must remain in $K$. It is easy to see  that condition \eqref{one sided} does not necessarily imply the uniqueness and follows immediately if, for instance, $f$ satisfies the local Lipschitz type condition of the form $\|f(x)-f(y)\|\leq\omega(\|x-y\|)$ for $x\in U\setminus K$ and $y\in K$. It may be shown that  the right semi-inner product $[\cdot,\cdot]_+$ in \eqref{one sided} can be replaced by the left semi-inner product $[\cdot,\cdot]_-$ and the assumption on  $\omega$ can be slightly relaxed.\\
\indent The viability of $K$ implies that $f$  is {\em tangent} to $K$, i.e. \begin{equation}\label{tangency}f(x)\in T_K(x)\;\;\text{for all}\;\; x\in\part K,\end{equation}
where $T_K(x)$ is the {\em contingent cone} (see \cite[Definition 4.1.1]{AF}).
In view of \cite[Theorem 3.5.7]{CNV} condition \eqref{tangency} along with \eqref{one sided} imply the invariance, too. In this case a solution starting at $x_0\in K$ is compared with an approximate solution living in $K$.\\
\indent If the set $K$ is {\em proximal} (see \cite[Def. 2.2.3]{CNV}), then in view of \cite{Volk1} (and comp. \cite[Theorem 4.3.1]{CNV}) conditions \eqref{tangency} and \eqref{one sided}, even in a slightly relaxed form (see \cite[Eq. (4.3.1)]{CNV}) imply the so-called {\em exterior tangency} condition saying that a {\em lower right Dini derivative} of the distance function $d_K=d(\cdot,K):=\inf_{k\in K}\|\cdot-k\|$ at $x$ in the direction $f(x)$ satisfies
\begin{equation}\label{e-tangency}D_+d_K\big(x;f(x)\big)\leq \omega\big(d_K(x)\big)\;\; \text{for}\;\;x\in U.\end{equation}
In \cite{Volk2} it is shown that the restrictive assumption concerning proximality of $K$ is  actually superfluous if $f$ is Lipschitz in the sense $\|f(x)-f(y)\|\leq\omega\big(\|x-y\|\big)$ for $x\in U\setminus K$, $y\in\part K$, where $\omega$ is a {\em uniqueness function}. Condition \eqref{e-tangency} actually entails the invariance (see \cite[Theorem 4.2.1]{CNV}).\\
\indent  The results described above represent two possible approaches to invariance. The first one relies on a controlled `monotonicity' (or a one-sided Lipschitz) condition \eqref{one sided}, which enables to compare an arbitrary solution to \eqref{equat+f} with another one surviving in $K$. The second approach does not use monotonicity. Instead the Lagrange type stability assumption \eqref{e-tangency} is employed. It requires that $d_K$ plays a role of a Lyapunov function of sorts that does not allow solutions to escape from $K$.\\
\indent  In the case of \eqref{equat-A+f} with a nontrivial $m$-accretive operator $A$ some important  results concerning invariance are also known. They can be found for instance in \cite{Pavel} or \cite{Shi} for linear $A$ and \cite{Bothe1996} for nonlinear $A$, see also \cite[Theorems 10.10.1 and 10.10.2]{CNV} and the references therein. Generally speaking, a closed set $K\subset \overline{D(A)}\cap\Omega$ is invariant with respect to \eqref{equat-A+f} if $f$ satisfies condition \eqref{one sided} above and either $K$ is viable or $f$ is $A$-tangent to $K$. See also \cite{cdp, cp2011} for results with set-valued $f$.\\
\indent It is our aim to deal with the `stability' approach rather, and  to get the analogue of the exterior tangency condition \eqref{e-tangency} in the situation of \eqref{equat-A+f}. In this context two recently published papers \cite{CaDaPraFran, CaDaPraFran2020} for linear $A$ have been direct inspirations for our considerations, where conditions outside $K$ expressed in terms of Dini derivatives in directions of $-A+f$ where considered. \\
\indent As it appears in applications, the constraint set $K$ is represented as a sublevel set of a certain locally Lipschitz potential $V:X\to\R$, i.e.
\begin{equation}\label{constr}
K=K_V:=\big\{x\in \overline{D(A)} \mid V(x)\leq 0\big\}
(\footnote{Note that the  representation \eqref{constr} does not restrict the generality since  is applicable for {\em any} closed set $K\subset \overline{D(A)}$ with $V=d_K$ or with $V=\Delta_K:=d_K-d_{X\setminus\wew K}$ if $\wew K\neq\emptyset$.}),
\end{equation}
and it seems natural to look for invariance conditions in terms of the constraining functional $V$.
The key ingredient of our approach relies on the analysis of behavior of $V$ along solutions of \eqref{equat-A+f}. The concept of the Dini {\em $A$-directional derivative} $D_AV (x;v)$ of $V$ at point $x\in\overline{D(A)}$ in the direction $v = f(x)$ is useful. We propose the following results without the reflexivity assumption on $X$  (see Sections \ref{sec2} and \ref{sec3} for the terminology and notation).
\begin{thm}\label{thm-inv1} If for every $z\in \partial K$ there are a neighborhood $U(z)\subset \Omega$ of $z$ and a uniqueness function $\omega$ such that
\begin{equation}\label{dense}
D_AV\big(x;f(x)\big)\leq\omega\big(V(x)\big)\;\; \text{for}\;\; x\in U(z)\setminus K,\end{equation} then $K$ is invariant with respect to \eqref{equat-A+f}.
\end{thm}
\noindent As the uniqueness function one often uses a linear function $\omega(t)=Ct$, where $C\geq 0$. In this case condition \eqref{dense} has the form
$$
D_A V(x;f(x)) \leq C V(x).
$$

\indent We also present a version of Theorem \ref{thm-inv1}, which is convenient in applications to partial differential equations of parabolic type (see 5.2 and 5.3).
\begin{thm}\label{thm-inv1-prime}
Suppose that for every $z\in \partial K$ there are a neighborhood $U(z)\subset \Omega$ of  $z$  and a uniqueness function $\omega$ such that, for any integral solution $u:[0,\tau)\to X$ of \eqref{equat-A+f} with $u(0)=z$ one has
\begin{equation}\label{dense-prime}
D_+(V\circ u)(t) \leq \omega \big( V\big(u (t)\big)\big), \ \  \mbox{ for a.e. } \ t\in (0,\tau), \ \mbox{ whenever  } u(t)\in U(z)\setminus K.
\end{equation}
Then $K$ is invariant with respect to \eqref{equat-A+f}.
\end{thm}
\noindent There are situations (e.g. for first order partial differential equations) when the verification of condition \eqref{dense} for $x$ from $\overline{D(A)}$ or condition \eqref{dense-prime} for any integral solution is not obvious if possible. For that reason, the following result seems to be more suitable.

\begin{thm} \label{thm-inv2}
Assume that $X$ is reflexive and:
\vspace{-3mm}
\begin{enumerate}
\item [\em(i)] for any $x\in\part K$, there is $\delta>0$  and a slow function $\beta$ such that $D(x,\delta)\subset\Omega$ and
\begin{equation}\label{slow war}\big[u-v,f(u)-f(v)\big]_+\leq\beta\big(\|u-v\|\big)\;\; \text{for}\;\; u,v\in D(x,\delta);\end{equation}
\item [\em (ii)] $f$ maps bounded sets into bounded ones;
\item [\em (iii)] for every $z\in \part K$ there is a neighborhood $U(z)\subset \Omega$ of $z$ such that
\begin{equation}\label{weak dense}D_AV \big(x;f(x)\big)\leq \omega\big(V(x)\big) \;\; \text{for}\;\;x\in D(A)\cap \big(U(z)\setminus K\big),\end{equation}
where $\omega$ is a nondecreasing uniqueness function.
\end{enumerate}
\vspace{-3mm}
Then $K$ is invariant with respect to \eqref{equat-A+f}.
\end{thm}
\noindent A notion of a slow function mentioned in the above theorem is provided in Subsection \ref{sf}. Note that condition \eqref{slow war} implies, in particular, that integral solutions to \eqref{equat-A+f} starting in a neighborhood of $\part K$ are locally unique. This is the price for relaxing condition \eqref{dense}.\\
\indent If $X$ is not reflexive then the following result will be useful when applying to partial differential equations (see Section 5.4).
\begin{thm}\label{thm-inv2-prime}
Assume that conditions {\em (i)} and {\em (ii)} from Theorem \ref{thm-inv2} are satisfied and that
\vspace{-3mm}
\begin{enumerate}
\item [\em (iii)'] for every $z\in \part K$ there are a neighborhood $U(z)\subset \Omega$ of $z$ and a nondecreasing uniqueness function $\omega$ such that, for any integral solution
$u:[0,\tau) \to X$ with $u(0)\in D(A)$, for a.e. $t\in (0,\tau)$ with $u(t)\in U(z)\setminus K$ one has
\begin{equation}\label{weak dense-prime}
D_+ (V\circ u) (t) \leq \omega\big(V\big(u(t) \big)\big).
\end{equation}
\end{enumerate}
\vspace{-3mm}
Then $K$ is invariant with respect to \eqref{equat-A+f}.
\end{thm}
\begin{rem} \label{how-to-compute}
{\em (1) Assume, additionally, that $X$ embeds continuously into  another Banach space $Y$ and that there is a quasi $m$-accretive operator $A_Y: D(A_Y)\to Y$ such that the part of $A_Y$ in $X$ is equal to $A$ and that $V$ can be extended to a differentiable $V_Y:Y\to \R$. Suppose also that, for any integral solution $u:[0,T]\to X$ of \eqref{equat-A+f},
$u\in W_{loc}^{1,1}((0,T],Y)$ (i.e. $u\in W^{1,1}([\delta,T],Y)$, for any $\delta\in (0,T)$), $u(t)\in D(A_Y)$ and
\begin{equation}\label{equat-A+f-Y}
\dot u(t)\in - A_Y u(t) + f(u(t)) \ \mbox{ for a.e. } \ t\in [0,T],
\end{equation}
i.e., in other words any integral solution to \eqref{equat-A+f} is a strong solution to \eqref{equat-A+f-Y} (see Subsection \ref{cp} (a) and \ref{cp-per}).
Then, for a.e. $t\in [0,T]$,
$$
D_+ (V\circ u)(t) = D_+ (V_Y\circ u)(t)=V_Y'(u(t))\dot u(t) =
V_Y'(u(t))\big(-v +f(u(t))\big)
$$
for some $v\in A_Y u(t)$. Therefore in order to verify \eqref{dense-prime} it is enough to show that, for all $x\in (U(z)\setminus K) \cap D(A_Y)$,
\begin{equation}\label{dens-prime-exp}
V_Y'(x)\big(-v+f(x)\big) \leq \omega \big(V(x)\big) \ \mbox{ for all } \ v\in A_{Y} x.
\end{equation}
\indent (2) Now assume that $X$ embeds continuously into a reflexive Banach space $Y$ and there exists  a quasi $m$-accretive operator $A_Y: D(A_Y)\to Y$ such that the part of $A_Y$ in $X$ is equal to $A$ and that $V$ can be extended to a differentiable functional $V_Y:Y\to \R$.
Moreover, suppose that $f$ is locally Lipschitz and bounded on bounded sets. Then one can prove (see Proposition \ref{reg-trick}) that if $u:[0,T]\to X$ is an integral solution with $u(0)\in D(A)$, then
$u\in W^{1,1}([0,T],Y)$, $u(t)\in D(A_Y)$ and
$$
\dot u(t)\in - A_Y u(t) + f\big(u(t)\big), \ \mbox{ for a.e. } \ t\in [0,T].
$$
In consequence, for a.e. $t\in [0,T]$,
$$
D_+ (V\circ u) (t) = D_+ (V_Y\circ u)(t) = V_Y'(u(t))\dot u(t)
= V_Y'(u(t)) \big(-v+f\big(u(t)\big) \big)
$$
for some $v\in A_Y u(t)$. Therefore in order to show that condition \eqref{weak dense-prime} holds it is sufficient to check \eqref{dens-prime-exp}.
}
\end{rem}
Results from Theorems \ref{thm-inv1} -- \ref{thm-inv2-prime} together with the above remarks  generalize results from \cite{CaDaPraFran} where linear operator $A$ was considered and $f$ was assumed to be globally quasi-dissipative. Here $A$ may be nonlinear and $f$ only continuous or locally quasi-dissipative. Moreover, we have come up with the invariance criteria that do not require the reflexivity of $X$.

The following results constitute analogs of the above theorems and refer to invariance of the interior of $K$. Our standing assumption is
\begin{equation}\label{constr^0}
K_V^0:=\big\{x\in \overline{D(A)} \mid V(x)< 0\big\}=\wew K_V\neq\emptyset.
\end{equation}
\begin{thm}\label{thm-inv-kv0-1} If for every $z\in \partial K$ there are a neighborhood $U(z)\subset \Omega$ of $z$ and a uniqueness function $\omega$ such that
\begin{equation}\label{dense1}D_AV\big(x;f(x)\big)\leq\omega\big(-V(x)\big)\;\; \text{for}\;\; x\in U(z)\cap K_V^0,\end{equation} then $K_V^0$ is invariant with respect to \eqref{equat-A+f}.
\end{thm}
\begin{thm} \label{thm-inv-kv0-2} Assume that $X$ is reflexive, assumptions (i) and (ii) from Theorem \ref{thm-inv2} hold and
\vspace{-3mm}
\begin{enumerate}
\item [\em (iv)] for every $z\in \part K$ there is a neighborhood $U(z)\subset \Omega$ of $z$ such that
\begin{equation}\label{weak dense1}D_AV\big(x;f(x)\big)\leq \omega\big(-V(x)\big) \;\; \text{for}\;\;x\in D(A)\cap (U(z)\cap K_V^0),\end{equation}
where $\omega$ is a nondecreasing uniqueness function.
\end{enumerate}
\vspace{-3mm}
Then $K_V^0$ is invariant with respect to \eqref{equat-A+f}.
\end{thm}
\indent Introducing additionally the {\em inwardness} condition \eqref{strict-inward} we obtain the following strict invariance results.
\begin{thm}\label{thm-strinv-1} If for every $z\in \partial K$ there are a neighborhood $U(z)\subset \Omega$ of $z$ and a uniqueness function $\omega$ such that
\begin{gather}\label{dense2}D_AV\big(x;f(x)\big)\leq\omega\big(|V(x)|\big)\;\; \text{for}\;\; x\in U(z),\\
\label{strict-inward}
D_AV\big(x;f(x)\big)<0\;\; \text{for}\;\; x\in U(z)\cap \part K,\end{gather}
then $K$ is strictly invariant with respect to \eqref{equat-A+f}.
\end{thm}
\begin{thm} \label{thm-strinv-2} Assume that $X$ is reflexive, assumptions (i) and (ii) of Theorem \ref{thm-inv2} hold,
\vspace{-3mm}
\begin{enumerate}
\item [\em (v)] for every $z\in \part K$ there is a neighborhood $U(z)\subset \Omega$ of $z$ such that
\begin{equation}\label{weak dense2}D_AV\big(x;f(x)\big)\leq \omega\big(|V(x)|\big) \;\; \text{for}\;\;x\in D(A)\cap U(z),\end{equation}
where $\omega$ is a nondecreasing uniqueness function, \end{enumerate}
and the strong inwardness condition \eqref{strict-inward} holds true.

\vspace{-3mm}
Then $K$ is strictly invariant with respect to \eqref{equat-A+f}.
\end{thm}
\noindent Theorems \ref{thm-inv1} -- \ref{thm-inv2-prime} will be proved in Section 3 and Theorems \ref{thm-inv-kv0-1} -- \ref{thm-strinv-2} in Section \ref{sec-str-inv}.\\
\indent Section \ref{applications} is devoted to applications. We discuss an evolution equation of the form \eqref{equat-A+f} under the presence of state dependent impulses -- see Subsection 5.1. In Subsection 5.2 we consider a fully nonlinear double obstacle problem
$$\left\{
\begin{array}{l}
u_t - \Delta_p u=f(x,u),\ \  x \in (0,l),\ t>0\\
u(0,t)= u(l,t) =0  \ \ t>0,\\
m(x)\leq u(x,t) \leq M(x),\  \  x\in (0,l), \, t>0
\end{array}
\right.
$$
where $\Delta_p$ is the $p$-Laplace operator and functions $m$ and $M$ represent two bodies, between which solutions are expected. Since the $p$-Laplace operator is nonlinear, results from \cite{CaDaPraFran} do not apply here. In Subsection 5.3 we study also the higher dimensional obstacle problem with the Laplace operator and in Section 5.4 we deal with invariance in the McKendrick age-structured population model. Let us note that, since the reflexivity of $X$ is not required in the applied abstract results (i.e. Theorems \ref{thm-inv1-prime} and \ref{thm-inv2-prime}), we are able to study partial differential equations as problems \eqref{equat-A+f} with $X$ being a space of continuous functions. Therefore we put neither growth restrictions nor global Lipschitz conditions on $f$.

\section{Preliminaries}\label{sec2}
 \noindent {\em Notation:} The notation used throughout the paper is standard. Given a metric space $(Y,d)$, $K\subset Y$ and $x\in Y$, $d_K(x):=\inf_{y\in K}d(x,y)$; by $\ov K$, $\wew K$ and $\part K$ we denote the closure, the interior and the boundary of $K$; $B(x,r)$ \big(resp. $D(x,r)$\big) is the open (resp. closed) ball around $x\in Y$ of radius $r>0$; $B(K,r)$ denotes the $r$-neighborhood of $K$, i.e. $B(K,r)=\{y\in Y\mid d_K(x)<r\}$. In what follows $\big(X,\|\cdot\|\big)$ is a real Banach space, $X^*$  stands for the  dual of $X$; $\la\cdot,\cdot\ra$ is the conjugation duality in $X$, i.e. if $x\in  X$, $p\in X^*$, then $\la p,x\ra:=p(x)$; by default $X^*$ is normed.
  The use of function spaces ($L^p$, Sobolev $W^{k,p}$, etc.), linear (unbounded in general) operators in Banach spaces, $C_0$ semigroups is standard. In particular, given real functions $u,v$, we put $u\vee v:=\max\{u,v\}$, $u\wedge v:=\min\{u,v\}$, $u_{\pm}=(\pm u)\vee 0$.

 We collect and recall some general concepts and relevant facts concerning evolution problems involving accretive operators.
 \subsection{General concepts}
 \subsubsection{Dini derivatives}\label{Dini} For any function $u:(a,b)\to\R$ one defines the {\em  Dini derivatives}
 $$D_\pm u(t):=\liminf_{h\to 0^\pm}\,\frac{u(t+h)-u(t)}{h}, \ \ \ D^\pm u(t):=\limsup_{h\to 0^\pm}\,\frac{u(t+h)-u(t)}{h}, \ \ \ \ \ \ t\in (a,b).$$
If  $f:\Omega\to\R$, where $\Omega\subset X$ is open, $x\in\Omega$ and $v\in X$, then the Dini directional derivatives at $x$ in the direction $v$ are given by
 $$D_\pm f(x;v):=\liminf_{h\to 0^\pm}\,\frac{f(x+hv)-f(x)}{h},\;\; \ \ D^\pm f(x;v):=\limsup_{h\to 0^\pm}\,\frac{f(x+hv)-f(x)}{h}.$$
 If the function $f$ is convex, then $D_- f(x;v)=D^- f(x;v)$ and $D_+ f(x;v)=D^+ f(x;v)$.
\subsubsection{Semi-inner products {\em (see e.g. \cite[Sec. 1.6]{CNV})}}\label{semi} If $x,y\in X$, then
we put
 $$[x,y]_\pm =\lim_{h\to 0^\pm}\frac{\|x+hy\|-\|x\|}{h},$$
 i.e. $[x,y]_\pm$ is the lower right (resp. left) Dini directional  derivative of $\|\cdot\|$ at $x$ in the direction of $y$ (see \cite[Ex. 1.6.1]{CNV} for properties of $[\cdot,\cdot]_\pm$).
\subsubsection{Uniqueness functions}\label{uf} A continuous function $\omega:[0,\infty)\to [0,\infty)$ such that $\omega(0)=0$ is a {\em uniqueness} (or a {\em Perron})  function if the only nonnegative solution to the problem $\dot u=\omega(u)$ on an interval $[0,\tau)$, $0<\tau\leq\infty$, such that $u(0)=0$ is the null function (some authors may consider different definitions, see, e.g. \cite{Walter}). Examples of uniqueness functions are numerous.

In the following lemma the first result is essentially due to Perron, the second one follows as an easy consequence of the first one, and the third one is proved for the sake of completeness. By a maximality of a solution $x$ we mean that $u(t)\leq x(t)$ for every other solution $u$ and every $t$ in a common interval of existence.
\begin{lm}\label{Perron} Let $\omega:[0,\infty)\to [0,\infty)$ be continuous and $\tau>0$.
 \vspace{-3mm}
 \begin{enumerate}
  \item[{\em (1)}] {\em (see, e.g. \cite[Th. I.4.10]{Coppel})} If $x$ is the maximal solution to the problem $\dot x=\omega(x)$ on the interval $[0,\tau]$, $u:[0,\tau]\to [0,\infty)$ is continuous, $u(0)\leq x(0)$ and  $Du\leq \omega(u)$  on $(0,\tau)$, where $D$ stands for any Dini derivative, then  $u(t)\leq x(t)$ for $t\in [0,\tau]$. In particular if $\omega$ is a uniqueness function and $u(0)=0$, then $u\equiv 0$ on $[0,\tau]$.
  \item[{\em (2)}] If $\omega$ is a uniqueness function, $u:[-\tau,0]\to (-\infty,0]$ is continuous, $u(0)=0$ and  $Du\leq \omega(-u)$  on $(-\tau,0)$, where $D$ stands for any Dini derivative,  then $u\equiv 0$ on $[-\tau,0]$.
  \item[{\em (3)}]  If $x$ is the maximal solution to the problem $\dot x=\omega(x)$ on the interval $[0,\tau]$, $u:[0,\tau]\to [0,\infty)$ is continuous and $u\in W_{loc}^{1,1}((0,\tau])$, $u(0)\leq x(0)$ and $\dot u(t)\leq \omega(u(t))$ for a.a $t\in (0,\tau]$, then $u\leq x$ on $[0,\tau]$. If $\omega$ is a uniqueness function, $u:[0,\tau]\to [0,\infty)$ is continuous and $u\in W^{1,1}_{loc}((0,\tau])$, $u(0)=0$ and  $\dot u\leq \omega(u)$ a.e. on $(0,\tau]$, then  $u\equiv 0$ on $[0,\tau]$.
      \end{enumerate}
      \vspace{-2mm}
\end{lm}

\begin{proof} We have the following property: {\em Let $\omega:[0,+\infty)\to [0,+\infty)$ be continuous and $\tau>0$. If $x$ is the maximal solution to the problem $\dot x=\omega(x)$ on $[0,\tau]$, $u:[0,\tau]\to [0,+\infty)$ is absolutely continuous, $\dot u\leq\omega(u)$ a.e. $[0,\tau]$ and $u(0)\leq x(0)$, then $u\leq x$ on $[0,\tau]$}. This, actually, follows from \cite[Theorem III.16.2]{Sz}.\\
\indent  Now let $t_0=\sup\{t\in [0,\tau]\mid u(s)\leq x(s)\; \text{for}\; 0\leq s\leq t\}$. Clearly $t_0\geq 0$. Suppose $t_0<\tau$. For $n\in\N$, let $y_n:[t_0,t_n]\to [0,+\infty)$, where $t_n>t_0$, be a solution to the problem $\dot y=\omega(y)+\frac{1}{n}$, $y(t_0)=x(t_0)+\frac{1}{n}$ on $[t_0,t_n]$.
Following \cite[Theorem I.4.9]{Coppel} or \cite[Lemma 1.8.1]{CNV}, we can assume that there is $\delta>0$ such that $t_0+\delta\leq t_n$ for all $n$. Moreover $y_n\geq y_{n+1}$ on $[t_0,t_0+\delta]$. Hence $y_n$ converges uniformly on $t_0,t_0+\delta]$ to a solution $y$ of the problem $\dot y=\omega(y)$, $y(t_0)=x(t_0)$. It is immediate to see that $y$ is the maximal solution; hence $y=x$ on $[t_0,t_0+\delta]$. On the other hand we have that $u(t_0)=x(t_0)<y_n(t_0)$ for any $n\in\N$. By continuity there is $t_0<t_0^n<t_0+\delta$ such that $u(t)\leq y_n(t)$ for $t_0\leq t\leq t_0^n$. Since $u\in W^{1,1}([t_0^n,\tau])$, $u$ is absolutely continuous on $[t_0^n,t_0+\delta]$ and
$\dot u\leq \omega(u)\leq \omega(u)+\frac{1}{n}$. Thus, by the above result, $u\leq y_n$ on $[t_0^n,t_0+\delta]$. Hence $u\leq y_n$ on the whole interval $[t_0,t_0+\delta]$, i.e. $u\leq x$ on this interval, too. This contradicts the choice of $t_0$ and proves the assertion.
\end{proof}

\subsubsection{Slow functions}\label{sf} A continuous nondecreasing function $\beta:[0,\infty)\to [0,\infty)$ is {\em slow} if there are $\eps>0$, $M>0$ and $\tau>0$ such that if $u:[0,\tau]\to [0,\infty)$ is continuous and
\begin{equation}\label{ineq-beta}u(t)\leq a+\int_0^t\beta(u(s))\,ds\;\;\text{for}\;\; t\in [0,\tau],\end{equation}
where $a\in [0,\eps]$, then one has $u(t)\leq a M$ for $t\in [0,\tau]$.

\indent If $\beta$ is a slow function, then the growth of a local solution of the problem $\dot u=\beta(u)$, $u(0)=a$, where $a$ is sufficiently small, is {\em slow} in the sense that it is `proportional' to the initial value $a$. This means that the use of slow functions provide arguments similar to those relying on the Gronwall inequality. It is clear that a slow function is a uniqueness function.\\
\indent The class of slow functions is large. For instance any continuous and nondecreasing function $\beta:[0,\infty)\to [0,\infty)$ such that $\beta(0)=0$, $\beta(x)>0$ if $x>0$ and there is $\eps>0$ such $\beta(\lambda x)\leq\lambda\beta(x)$ for all $0\leq\lambda\leq\eps$ and $x\geq 0$, is slow. In particular any function $\beta$ of the form $\beta(x)=x\gamma(x)$, $x\geq 0$, where $\gamma:[0,\infty)\to [0,\infty)$ is continuous nondecreasing and $\gamma(x)>0$ when $x>0$, is slow. For instance slow are  the following functions $\beta(x)=cx^\alpha$, where $c\geq 0$ and  $\alpha\geq 1$, or $\beta(x)=x\ln(1+x)$ for $x\geq 0$.\\
\indent All these examples are subsumed by the following general result.
\begin{prop}\label{slow}
If $\beta:[0,+\infty)\to [0,+\infty)$ is a continuous and nondecreasing function such that
\begin{equation}\label{liminf-condition}
\liminf_{x\to 0^+} \frac{x}{\beta(x)} >0,
\end{equation}
then $\beta$ is slow.
\end{prop}
\begin{proof}
For any $a>0$, define $B_a:[1,+\infty) \to [0,B_{a,\infty})$
by
$$
B_a (r):= \int_{1}^{r} \frac{a}{\beta(as)} \, ds, \ r>1,
$$
where $B_{a,\infty}:= \lim_{r\to +\infty} B_a(r)$.  Take any $\Gamma$ such that
$$
\liminf_{x\to 0^+} x/\beta(x) >\Gamma>0.
$$
In view of \eqref{liminf-condition}, there exists $\varepsilon>0$ such that
$$
x / \beta (x) > \Gamma \ \mbox{ for all } \ x\in (0,2\varepsilon].
$$
Therefore, if $s \in (1,2)$ and $a\in (0,\varepsilon]$, then
$$
\frac{a}{\beta(as)} \geq \frac{\Gamma}{s},
$$
which implies
$$
B_a (2) \geq  \int_{1}^{2} \Gamma s^{-1}\, ds \geq \Gamma/2 \  \mbox{ for all }\  a\in (0,\varepsilon],
$$
that is, with $M:=2$ and $\tau:=\Gamma/2$,
\begin{equation}\label{key-ineq}
B_a(M)\geq \tau \ \mbox{ for all } a\in (0,\varepsilon].
\end{equation}
Assume that $a\in (0,\varepsilon]$ and $u:[0,\tau]\to [0,+\infty)$ satisfies \eqref{ineq-beta}.
If we define $x_a:[0,\tau]\to [0,+\infty)$ by $x_a(t):=u(t)/a$ then
\begin{equation}\label{assump-restated}
x_a (t) \leq y_a(t):= 1+\int_{0}^{t} a^{-1}\beta(ax_a(s))\, d s.
\end{equation}
By definition and the monotonicity of $\beta$, we get
$$
a\dot y_a(t) = \beta(ax_a(t)) \leq \beta(ay_a(t)) \ \mbox{ for all } t\in [0,\tau].
$$
This yields
$$
\int_{0}^{t} \frac{a\dot y_a(s)}{\beta(ay_a(s))}\, ds \leq t  \ \mbox{ for all } t\in [0,\tau],
$$
which, by the change of variables and \eqref{key-ineq}, implies
$$
B_a(y_a(t))=\int_{1}^{y_a(t)} \frac{a}{a\beta(s)} \, ds = \int_{0}^{t} \frac{a\dot y_a(s)}{\beta(ay_a(s))}\, ds\leq t\leq \tau \leq B_a (M) \ \mbox{ for all } \ t\in [0,\tau].
$$
Hence, by the monotonicity of $B_a$, we get
$$
y_a(t)\leq M \ \mbox{ for all } \ t\in [0,\tau).
$$
Finally, in view of \eqref{assump-restated}, we see that, for all $a\in (0,\varepsilon]$,
$$
u(t)/a=x_a(t) \leq y_a(t) \leq M \ \mbox{ for all } \ t\in [0,\tau],
$$
which ends the proof.
\end{proof}
\begin{rem}\label{slow rem} {\em Observe that if a function $\beta:[0,\infty)\to [0,\infty)$ is slow, then for any $a,b\geq 0$, the function $[0,\infty)\ni x\mapsto a\beta(x)+bx$ is slow, too.}
\end{rem}

 \subsection{Nonlinear evolution associated with accretive operators}
  \subsubsection{Accretive operators}\label{opac} Let $A:D(A)\multimap X$, where $D(A)\subset X$ and  $X$ is a Banach space, be a set-valued operator, i.e. $\emptyset\neq Ax\subset X$ for $x\in D(A)$. Let $\Gr(A):=\{(x,u)\in X\times X\mid x\in D(A),\;u\in Ax\}$ be the graph of $A$.\\
\indent (a) $A$ is {\em accretive} if $[x-y,u-v]_+\geq 0$ for all $(x,u), (y,v)\in\Gr(A)$. $A$ is {\em $m$-accretive} if it is accretive and $\mathrm{Range}(I+\lambda A):=\big\{y\in X\mid y\in x+\lambda Ax\; \text{for some}\; x\in D(A)\big\}=X$ for some (equivalently for all) $\lambda>0$.\\
\indent (b) $A$ is {\em $\alpha$-accretive} (resp. $\alpha$-$m$-accretive), $\alpha\in\R$, if $\alpha I+A$ is accretive (resp. $m$-accretive);  hence $A$ is $\alpha$-accretive if and only if $[x-y,u-v]_+\geq -\alpha\|x-y\|$ for all $(x,u), (y,v)\in\Gr(A)$. $A$ is {\em quasi $m$-accretive} if it is $\alpha$-$m$-accretive for some $\alpha\in\R$.\\
\indent (c) If  $A$ is quasi $m$-accretive, then $\Gr(A)$ is closed and, in particular, the set $Ax$, $x\in D(A)$, is closed. If the dual $X^*$ is uniformly convex, then $Ax$ is convex. If both $X$ and $X^*$ are uniformly convex, then the closure $\overline{D(A)}$ is convex and for each $x\in D(A)$ and $w\in X$ there is a unique element $(Ax-w)^0\in Ax-w$ of minimal norm, i.e. $\|(Ax-w)^0\|=|Ax-w|:=\inf_{u\in Ax}\|u-w\|$.\\
\indent (d) In view of the Lumer theorem a  {\em linear} operator $A:D(A)\to X$ is $\alpha$-$m$-accretive if and only if $-A$ is a closed densely defined generator of a strongly continuous semigroup of linear operators $\{e^{-tA}\}_{t\geq 0}$ such that $\|e^{-tA}\|\leq e^{t\alpha}$ for $t\geq 0$.\\
\indent (e) If $A$ is $\alpha$-$m$-accretive, $\lambda>0$ with $\lambda\alpha<1$, then the {\em resolvent}
$J_\lambda=J_\lambda^A:=(I+\lambda A)^{-1}:X\to D(A)$ and the {\em Yosida approximation} $A_\lambda=\lambda^{-1}(I-J_\lambda):X\to X$  are well-defined, single-valued, and
\begin{equation}\label{resol}
\begin{split}&\|J_\lambda x-J_\lambda y\|\leq (1-\lambda\alpha)^{-1}\|x-y\|,\;\;A_\lambda x\in AJ_\lambda x\;\; \text{for all}\;\; x,y\in X,\\
&\lim_{\lambda\to 0^+}J_\lambda x=x\;\; \text{for}\;\; x\in\overline{D(A)}.
\end{split}
\end{equation}
\subsubsection{Cauchy problems}\label{cp} Assume $A$ is an $\alpha$-$m$-accretive operator.  Let $T>0$, $w\in L^1\big([0,T],X\big)$ and consider the problem
\begin{equation}\label{cauchy-accretive-w}
\left\{ \begin{array}{l}
\dot u(t) \in - Au(t)+w(t),\;  t\in [0,T],\\
u(0)=x\in \overline{D(A)}.
\end{array} \right.
\end{equation}
\indent (a) A continuous function $u:[0,T]\to X$ is a  {\em strong} solution to \eqref{cauchy-accretive-w} if  $u\in W^{1,1}_{loc}\big((0,T],X\big)$, $u(t)\in D(A)$, $u(0)=x$ and $\dot u(t)-w(t)\in -Au(t)$ for a.a. $t\in (0,T]$; here $\dot u(t)$ stands for the ordinary strong derivative; the formula makes sense since  $u$ is differentiable a.a.\\
\indent (b) If $X$ is reflexive, $x\in D(A)$ and $w\in W^{1,1}\big([0,T],X\big)$, then \eqref{cauchy-accretive-w} has a unique strong solution $u\in W^{1,\infty}\big([0,T],X\big)$; if both spaces $X$ and $X^*$ are uniformly convex and $x\in D(A)$, then the strong solution $u$ is right differentiable, $\dot u_+$ is right continuous and $\dot u_+(t)+\big(Au(t)-w(t)\big)^0=0$ for a.a. $t\in (0,T)$; if in addition $w$ is continuous, then  $\dot u_+(0)=\big(-Ax+w(0)\big)^0$  (see \cite[Theorems 4.5, 4.6]{Barbu}).\\
\indent (c) \cite[Remarks to Corollary 4.2]{Barbu} A continuous function $u:[0,T]\to X$  is an {\em integral} solution to \eqref{cauchy-accretive-w} if $u(0)=x$ and for any $0\leq s\leq t\leq T$ and $(y,v)\in\Gr(A)$,
\begin{equation}\label{def-integral}e^{-t\alpha}\big\|u(t)-y\big\|\leq e^{-s\alpha}\big\|u(s)-y\big\|+\int_s^te^{-z\alpha}\big[u(z)-y,w(z)-v\big]_+\,dz.\end{equation}
It is known that a strong solution is an integral one. By \cite[Theorem 4.2, Corollary 4.2]{Barbu} \eqref{cauchy-accretive-w} admits a {\em unique} integral solution denoted by $u=u_A(\cdot
;x,w):[0,T]\to X$ (or $u(\cdot;x,w)$ if $A$ is the default operator) and $u(t)\in\overline{D(A)}$. \\
\indent (d) Given $x_1,x_2\in \overline{D(A)}$, $w_1, w_2\in L^1\big([0,T],X\big)$ the {\em Benil\'an inequality} holds: for $0\leq s\leq t\leq T$
\begin{align}\label{semi-equality}e^{-t\alpha}\big\|u_1(t)-u_2(t)\big\|&\leq e^{-s\alpha}\big\|u_1(s)-u_2(s)\big\|+\int_s^te^{-z\alpha}\big[u_1(z)-u_2(z),w_1(z)-w_2(z)\big]_+
\,dz\\&\leq e^{-s\alpha}\big\|u_1(s)-u_2(s)\big\|+\int_s^te^{-z\alpha}\big\|w_1(z)-w_2(z)\big\|\,dz,\nonumber\end{align}
where $u_i:=u_A (\cdot;x_i,w_i)$, $i=1,2$, (see \cite[eq. (4.14)]{Barbu}.\\
\indent (e) If $w\equiv 0$, then \eqref{cauchy-accretive-w} has a unique integral solution $u_A (\cdot;x,0)$ defined on $[0,\infty)$ (i.e. $u_A(\cdot;x,0)$ is an integral solution on $[0,T]$ for any $T>0$) and  the following {\em Crandall-Liggett formula}, see \cite[Theorem 4.3]{Barbu}, holds: for any $x\in X$ and $t\geq 0$,
\begin{equation*}
u_A (t;x,0)=\lim_{n\to\infty}J_{t/n}^nx\;\; (\footnote{Here $J_0(x)=x$, $x\in X$.}).
\end{equation*}
Let us put
$$
S_A(t)x:=u_A (t;x,0),\;\; t\geq 0.
$$
Then, for any $t\geq 0$, $S_A(t):\overline{D(A)}\to \overline{D(A)}$, $\big\|S_A(t)x-S_A(t)y\big\|\leq e^{t\alpha}\|x-y\|$ or all $x,y\in\overline{D(A)}$. The family $\{S_A(t)\}_{t\geq 0}$ is a (strongly continuous) semigroup of continuous maps, i.e. for any $x\in\overline{D(A)}$ the map $[0,\infty)\ni t\mapsto S_A(t)x$ is continuous,
$S_A(0)=I$ on $\overline{D(A)}$ and $S_A(t+s)=S(t)\circ S_A(s)$ for any $t,s\geq 0$ (see \cite[Proposition 4.2]{Barbu}).
\begin{rem}\label{dod info} {\em Let us derive some immediate consequences of the above facts.\\
\indent (a) If $x\in D(A)$, then $S_A(\cdot)x:[0,\infty)\to X$ is Lipschitz continuous on every compact interval $[0,\tau]$, $T>0$. To see this fix $\tau>0$ and take $0\leq s<t\leq \tau$; then, in view of \eqref{semi-equality}, \eqref{def-integral}, for any $v\in Ax$
\begin{align*}\big\|S_A(t)x-S_A(s)x\big\|&=\big\|S_A(s)S_A(t-s)x-S_A(s)x\big\|\leq e^{s\alpha}\big\|S_A(t-s)x-x\big\|\\
&\leq e^{s\alpha }\int_0^{t-s}e^{(t-s-z)\alpha}\|v\|\,dz\leq |Ax|\int_0^{t-s}e^{(t-z)\alpha}\,dz\leq \ell(t-s)\end{align*}
for some $\ell>0$. The same is true for integral solutions $u_A(\cdot,x,w)$ with {\em constant} $w$.\\
\indent  (b) If $X$ is reflexive and $x\in D(A)$, then $u=S_A(\cdot)x\in W^{1,\infty}_{loc}\big([0,\infty),X\big)$, $u(t)\in D(A)$  and $\dot u(t)\in-Au(t)$ for a.a. $t\geq 0$.\\
\indent (c) Suppose that $A$ is linear. Then $S_A(t)=e^{-tA}$ for any $t\geq 0$. Moreover
$u=u_A (\cdot;x,w)$ is a mild solution of \eqref{cauchy-accretive-w} in the sense of the Duhamel formula, i.e.
\begin{equation}\label{duhamel}
u(t)=e^{-tA}x+\int_0^te^{-(t-s)A}w(s)\,ds,\; t\geq 0.
\end{equation}
\indent (d) Let $u=u_A (\cdot;x,w)$ and fix a small $h>0$.
For all $0\leq s\leq t\leq T-h$ and $(z,v)\in\Gr(A)$, one has
\begin{align*}\big\|u(t+h)-z\big\|&\leq e^{(t-s)\alpha}\big\|u(s+h)-z\big\|+\int_{s+h}^{t+h}e^{(t+h-\xi)\alpha}\big[u(\xi)-z,w(\xi)-v\big]_
+\,d\xi\\
&=e^{(t-s)\alpha}\big\|u(s+h)-z\big\|+\int_s^te^{(t-\xi)\alpha}\big[u(\xi+h)-z,w(\xi+h)-v\big]_
+\,d\xi.\end{align*}
Let $y(t):=u(t+h)$, $t\in [0,T-h]$. We have
$$\big\|y(t)-z\big\|\leq e^{(t-s)\alpha}\big\|y(s)-z\big\|+\int_s^te^{(t-\xi)\alpha}\big[y(\xi)-z,\bar w(\xi)-v\big]_+\,d\xi,$$
where $\bar w(t):=w(t+h)$, $t\in [0,T-h]$. Thus $y$ is an integral solution to the problem $\dot y\in-Ay+\bar w$, $y(0)=u(h)$, i.e. $y=u_A \big(\cdot;u(h),\bar w\big)$.  In other words we have the formula
\begin{equation}\label{wzorek}
u_A (\cdot+h;x,w)=u_A \big(\cdot;u_A(h;x,w),w(\cdot+h)\big),\;t\in [0,T-h]\end{equation}
and, in view of \eqref{semi-equality}, for any $0\leq s\leq t\leq
T-h$
\begin{align}\nonumber\|u(t+h)-u(t)\|&\leq e^{(t-s)\alpha}\big\|u(s+h)-u(s)\big\|+\int_s^te^{(t-z)\alpha}\big[u(z+h)-u(z),w(z+h)-w(z)\big]_
+\,dz\\
&\leq e^{(t-s)\alpha}\big\|u(s+h)-u(s)\big\|+\int_s^te^{(t-z)\alpha}\big\|w(z+h)-w(z)\big\|\,dz.
\end{align}
}\end{rem}

\subsubsection{Continuous perturbations} \label{cp-per}
Let $A$ be an $\alpha$-$m$-accretive operator and $f:\Omega\to X$, where $\Omega\subset X$, be continuous. A continuous $u:[0,T]\to\Omega$, where $T>0$, is an {\em integral} (resp. {\em strong}) solution to \eqref{equat-A+f}, i.e.
\begin{align}\label{cauchy-accretive-f}
&\left\{ \begin{array}{l}
\dot u  \in -Au+f(u)\\
u(0)=x\in \overline{D(A)}\cap\Omega,
\end{array} \right.
\end{align}
if $u$ is an {\em integral} (resp. {\em strong}) {\em solution}  to \eqref{cauchy-accretive-w} with $w:=f\circ u$ (\footnote{This definition makes sense since here $w\in L^1([0,T],X)$. More generally (if $f$ is not assumed to be continuous), we can say that $u$ is an integral solution of \eqref{cauchy-accretive-f} if $f\circ u\in L^1([0,T],X)$ and $u$ is an integral solution to \eqref{cauchy-accretive-w} with $w=f\circ u$.
}). A continuous function $u:[0,\tau)\to X$, $0<\tau\leq\infty$, is an integral solution to \eqref{cauchy-accretive-f} if
for any $0<T<\tau$, $u$ restricted to $[0,T]$ is an integral solution to \eqref{cauchy-accretive-f} on $[0,T]$. An integral solution $u:[0,\tau)\to X$ is {\em noncontinuable} if it has no extension to a solution defined on the interval $[0,\tau')$ with $\tau'>\tau$.
\begin{rem}\label{semigroup prop}{\em (1) Along with \eqref{cauchy-accretive-f} consider the problem
\begin{align}\label{cauchy-B}
&\left\{ \begin{array}{l}
\dot u  \in -Bu+g(u)\\
u(0)=x\in \overline{D(A)},
\end{array} \right.
\end{align}
where $B=A+\alpha I$ and $g(u)=\alpha u+f(u)$, $u\in\Omega$. $B$ is $m$-accretive. We claim that integral solutions to \eqref{cauchy-accretive-f} and \eqref{cauchy-B} coincide. Clearly it is sufficient to show that
if $u:[0,T]\to\Omega$ solves \eqref{cauchy-B}, then it solves \eqref{cauchy-accretive-f}. By definition $u$ is an integral solution to $\dot u\in -Bu+v$, $u(0)=x$, where $v=g\circ u$. In view of \cite[Chapter 4.1]{Barbu} $u$ is the uniform limit of $\eps$-DS-approximate solutions $u^\eps$ as $\eps\to 0^+$, i.e. $\sup_{t\in [0,T]}\big\|u(t)-u^\eps(t)\big\|<\eps$. Here by an
$\eps$-DS-approximate solution to \eqref{cauchy-B} we mean a step function $u^\eps$ with
$$u^\eps(0)=x,\,u^\eps(t)=u_k\;\; \text{on}\;\; (t_{k-1},t_k\wedge T]\;\; \text{for}\;\; k=1,...,m,$$
where $0=t_0<t_1<...<t_{m-1}<T\leq t_m$, $t_k-t_{k-1}\leq \eps$ for $k=1,...,m$, and
$$\frac{u_k-u_{k-1}}{t_k-t_{k-1}}\in -Bu_k+v_k\;\; \text{for}\;\; k=1,...,m,$$
with $v_1,...,v_m\in X$ such that
$$\sum_{k=1}^m\int_{t_{k-1}}^{t_k\wedge T}\big\|v_k-v(t)\big\|\leq\eps.$$
It is clear that for all $k=1,...,m$,
$$\frac{u_k-u_{k-1}}{t_k-t_{k-1}}\in -Au_k+v_k-\alpha u_k,$$
and  for $t\in [t_{k-1},t_k\wedge T]$,
$$\big\|(v_k-\alpha u_k)-w(t)\big\|\leq \big\|v_k-v(t)\big\|+\big\|\alpha u(t)-\alpha u_k\big\|\leq \big\|v_k-v(t)\big\|+|\alpha|\eps.$$
Therefore
$$\sum_{k=1}^m\int_{t_{k-1}}^{t_k\wedge T}\big\|v_k-\alpha u_k-w(t)\big\|\leq\eps\big(1+|\alpha|T\big).$$
This shows that $u^\eps$ is an $\eps\big(1+|\alpha|T\big)$-DS-approximate solution of \eqref{cauchy-accretive-w} and, hence, $u$ is an integral solution to \eqref{cauchy-accretive-f}. As a conclusion we see while studying \eqref{cauchy-accretive-f} that one can shift the part $\alpha I$ from the perturbation term to the accretive operator and vice-versa. In particular we may, without any loss of generality, consider {\em only} $m$-accretive operators $A$.\\
\indent (2) Suppose $u:[0,\tau)\to\Omega$, where $\tau>0$, is an integral solution of \eqref{cauchy-accretive-f} and there is $M>0$ such that $\big\|f(u(t))\big\|\leq M$ for any $0\leq t<\tau$. Putting $w(t):=f\big(u(t)\big)$, $t\in [0,\tau)$, we have $w\in L^1\big([0,\tau],X\big)$. Hence there is a continuous function $\bar u:[0,\tau]\to X$ being an integral solution on $[0,\tau]$ of the problem $\dot u\in-Au+w$, $u(0)=x$. Evidently
$u=\bar u$ on $[0,\tau)$. Hence $\lim_{t\to \tau^-}u(t)=\bar u(\tau)\in\overline\Omega$ exists.\\
\indent (3) Let $u:[0,T]\to X$ be a solution to \eqref{cauchy-accretive-f}, i.e. $u=u_A (\cdot;x,w)$, where $w=f\circ u$. Let $0<h<T$. Since
$w(\cdot+h)=f\circ u(\cdot+h)$ on $[0,T-h]$, we see in view of \eqref{wzorek} that $y=u(\cdot+h):[0,T-h]\to X$ is a solution  to the problem $\dot y\in-Ay+f(y)$, $y(0)=u(h)$. This proves a  `semigroup' property of sorts of integral solutions to \eqref{cauchy-accretive-f}. Namely, putting
$$S_A(t;f)(x):=\{u(t)\mid u\; \text{is an integral solution to \eqref{cauchy-accretive-f}}\},\;\; t\in [0,T)$$
we have
$$S_A(t+h;f)(x)=S_A(t;f)\big(S_A(h;f)(x)\big)$$
 for $t,h\in [0,T]$ such that $t+h\leq T$.\hfill $\square$
}\end{rem}

\indent The following result seems to be a well-known folklore. It seems, however, to be more convenient than the corresponding result \cite[Theorem 4.8]{Barbu}. We sketch a proof for the sake of completeness.

\begin{prop}\label{prop lip} Assume that $f:\Omega\to X$ is locally Lipschitz continuous, $A:D(A)\multimap X$ is $m$-accretive and $x\in \overline{D(A)}\cap\Omega$. Then:
\begin{enumerate}
\item [{\em (a)}] there  is $T>0$ and a unique integral solution $u:[0,T]\to \Omega$ of \eqref{cauchy-accretive-f}.
\item [{\em (b)}] If $u:[0,T]\to\Omega$ is an integral solution of \eqref{cauchy-accretive-f} with  $x\in D(A)\cap\Omega$, then $u$ is Lipschitz continuous.
\item [{\em (c)}] If $X$ is reflexive, $u:[0,T]\to\Omega$ is an integral solution of \eqref{cauchy-accretive-f} with $x\in D(A)\cap\Omega$, then $u$ is a strong solution and
$u\in W^{1,\infty}\big([0,T],X\big)$.
\end{enumerate}
\end{prop}
\begin{proof} (a) There is $R>0$ such that $D:=D(x,R)\subset \Omega$ and $f$ is Lipschitz with the Lipschitz constant $\ell>0$ on $D$. Take $y\in D(A)\cap D(x,R/3)$ and $p\in Ay$. Let $M:=\sup_{u\in D}\big\|f(u)-p\big\|$ and
$$Y:=\big\{u\in C\big([0,T],X\big)\mid u(t)\in D,\;t\in [0,T]\big\},$$
where $T=\frac{R}{3M}$. Let $Y$ be endowed with the  complete metric
$$d(u,v):=\sup_{t\in [0,T]}e^{-\ell t}\big\|u(t)-v(t)\big\|,\;\; u,v\in Y.$$ Consider a map $N:Y\to C\big([0,T],X\big)$ given by
$$[Nu](t)=u_A (t;x,w_u),\; u\in Y,$$
where $w_u:=f\circ u$, i.e. $Nu$ is the integral solution to \eqref{cauchy-accretive-f} with $w=w_u$. This map is well-defined since $w_u\in C\big([0,T],X\big)\subset L^1\big([0,T],X\big)$. Actually
$N:Y\to Y$ and it is a (Banach) contraction. Indeed: for $u\in Y$ and $t\in [0,T]$
$$\big\|w_u(\tau)-p\big\|=\big\|f(u(\tau))-p\big\|\leq M,\;\; \tau\in [0,t],$$
and, in view of \eqref{def-integral}
\begin{align*}
\big\|[Nu](t)-x\big\|&\leq\big\|[Nu](t)-y\big\|+\|x-y\|\leq 2\|x-y\|+\int_0^t\big\|w_u(z)-p\big\|\,dz\leq
\frac{2}{3}R+MT\leq R,
\end{align*}
i.e. $Nu\in Y$. For $u,v\in Y$ in view of \eqref{semi-equality} we have
\begin{align*}\big\|[Nu](t)-[Nv](t)\big\|&\leq \int_0^t\big\|w_u(z)-w_v(z)\big\|\,dz\leq
\ell\int_0^t\big\|u(z)-v(z)\big\|\,dz\\
&\leq \ell d(u,v)\int_0^te^{\ell z}\,dz=(e^{\ell t}-1)d(u,v)\end{align*}
and thus $d(Nu,Nv)\leq cd(u,v)$ with $c=1-e^{-\ell T}<1$. Hence there is $u\in Y$ such that $Nu=u$, i.e. $u$ is an integral solution to \eqref{cauchy-accretive-f}. Its uniqueness is straightforward.\\
\indent (b) To show the Lipschitz continuity of the solution $u=u_A (\cdot;x,w):[0,T]\to \Omega$, where $x\in D(A)$ and $w:=f\circ u$, take $0\leq t<s\leq T$ and let $h:=s-t$. Arguing as in Remark \ref{dod info} (d) we have by \cref{def-integral,semi-equality,wzorek}  that for $v\in Ax$
\begin{align*}\big\|u(s)-u(t)\big\|&=\big\|u(t+h)-u(t)\big\|\leq \int_0^h\big\|w(z)-v\big\|\,dz+\int_0^t
\Big\|f\big(u(z+h)\big)-f\big(u(z)\big)\Big\|\,dz\\
&\leq ch+\ell \int_0^t \big\|u(z+h)-u(z)\big\|\,dz\end{align*}
where $c$ is a suitably chosen constant  and $\ell$ is the Lipschitz constant of $f$ restricted to the compact set $u\big([0,T]\big)$.  The Gronwall inequality implies that
\[\|u(s)-u(t)\|\leq ce^{LT}(s-t).\]
\indent (c) Observe that $w=f\circ u:[0,T]\to X$ is absolutely continuous as the superposition of two locally Lipschitz functions. This implies that $w\in W^{1,1}([0,T],X)$ in view of the Komura theorem (see \cite{Komura}). The assertion follows from Subsection \ref{cp} (b).
\end{proof}

\begin{thm}\label{main existence} Let $A:D(A)\multimap X$ be $m$-accretive and suppose that $f:\Omega\to X$ is continuous and
\begin{equation}\label{slow one sided}\big[u-v,f(u)-f(v)\big]_+\leq\beta\big(\|u-v\|\big)\;\; \text{for}\;\; u,v\in\Omega,\end{equation}
where $\beta$ is a slow function. Moreover assume that $f$ maps bounded sets into bounded ones. For any $x_0\in\overline{D(A)}\cap\Omega$ and $0<r<R$ such that $D(x_0,R)\subset\Omega$, there is $T>0$ such that for any $x\in D(x_0,r)$ the problem \eqref{cauchy-accretive-f} has a unique solution on $[0,T]$. If $x\in D(x_0,r)\cap D(A)$, then this solution is Lipschitz continuous.
\end{thm}
\begin{proof}  In view of \cite[Lemma 1]{Lasota} for any $n\in\N$ there is a locally Lipschitz $f_n:\Omega\to X$ such that $\big\|f(x)-f_n(x)\big\|\leq \frac{1}{2n}$ for $x\in\Omega$. Then for any $n\in\N$ and $u,v\in\Omega$
\begin{equation}\label{1}\big[u-v,f_n(u)-f_n(v)\big]_+\leq\beta\big(\|u-v\|\big)+\frac{1}{n}\leq \beta\big(\|u-v\|\big)+1. \end{equation}
Take $x_0\in\overline{D(A)}\cap\Omega$, $0<r<R$ such that $D(x_0,R)\subset\Omega$. Let $\eps<R-r$. Since $J_\lambda x_0\to x_0$ as $\lambda\to 0^+$ we find $\lambda_0>0$ such that $\|x_0-y_0\|<\eps$, where $y_0=J_{\lambda_0}x_0$. Let $v_0:=A_{\lambda_0}x_0$; then $y_0\in\Omega$ and $v_0\in Ay_0$ \big(see Subsection \ref{opac} (e)\big). Let $\gamma:=\big\|v_0-f(y_0)\big\|+2$ and $r_0=r+\eps$.\\
\indent Let us consider the (scalar) problem
\begin{equation}\label{sca}\dot z=\beta(z)+\gamma,\; z(0)=r_0,\end{equation}
and let $z:[0,\tau)\to [0,+\infty)$, where $0<\tau\leq\infty$, be the maximal solution to \eqref{sca}. There exists $0<T< \tau$ ($\tau$ comes from the definition of a slow function, see Subsection \ref{sf}) such that $z(t)\leq R$ for all $t\in [0,T]$.\\
\indent Take any $n\in\N$, an arbitrary $x\in D(x_0,r)$ and a noncontinuable integral solution $u_n$ to the problem
\begin{align}\label{cauchy-accretive-fn}
&\left\{ \begin{array}{l}
\dot u  \in -Au+f_n(u)\\
u(0)=x\in \overline{D(A)},
\end{array} \right.
\end{align}
defined on $[0,\tau_n)$.\\
\indent We claim that $T\leq \tau_n$ for any $n\in\N$. Let
$$g(t):=\big\|u_n(t)-y_0\big\|,\; 0\leq t<\tau_n.$$
Then, by definition, for  $0\leq t<s<\tau_n$,
\begin{align*}g(s)=& \big\|u_n(s)-y_0\big\|\leq \big\|u_n(t)-y_0\big\|+\int_t^s\big[u_n(z)-y_0,f_n\big(u_n(z)\big)-v_0\big]_+\,dz\\
=& g(t)+\int_t^s\big[u_n(z)-y_0,f_n\big(u_n(z)\big)-f_n(y_0)+f_n(y_0)-v_0\big]_+\,dz\\
\leq& g(t)+\int_t^s\left(\beta\big(\|u_n(z)-y_0\|\big)+1+\big\|f_n(y_0)-v_0\big\|\right)\,dz\\
\leq& g(t)+\int_t^s\left(\beta\big(g(z)\big)+1+\big\|f(y_0)-v_0\big\|+\big\|f_n(y_0)-f(y_0)\big\|\right)
\,dz\\
\leq& g(t)+\int_t^s\left(\beta\big(g(z)\big)+\gamma\right)\,dz.
\end{align*}
Hence
$$D^+g(t)\leq \beta\big(g(t)\big)+\gamma,\; 0\leq t<\tau_n.
$$
Since $g(0)=\big\|u_n(0)-y_0\big\|=\|x-y_0\|\leq\|x_0-y_0\|+\|x_0-x\|\leq r_0$ this, in view of Lemma \ref{Perron}, implies that
$$g(t)\leq z(t)\;\; \text{for}\;\; 0\leq t<\tau_n\wedge \tau.$$
If, for some $n\in\N$,  $\tau_n<T$, then the set $\big\{u_n(t)\mid 0\leq t<\tau_n\big\}$ is bounded and, hence,  so is the set $\big\{f\big(u_n(t)\big)\mid 0\leq t<\tau_n\big\}$ as well as the set $\big\{f_n\big(u_n(t)\big)\mid 0\leq t<\tau_n\big\}$. This, in view of Remark \ref{semigroup prop} (2) implies that $\lim_{t\to\tau_n}u_n(t)$ exists., i.e. $u_n$ is continuable. The contradiction proves our claim.\\
\indent Now we shall show that the functional sequence $(u_n)$ converges uniformly. In view of \eqref{semi-equality} we have that for any $0\leq t\leq T$
\begin{align*}\big\|u_n(t)-u_m(t)\big\|&\leq\int_0^t\Big[u_n(z)-u_m(z),f_n\big(u_n(z)\big)
-f_m\big(u_m(z)\big)\Big]_+\,dz\\
&\leq \int_0^t\left(\Big[u_n(z)-u_m(z),f\big(u_n(z)\big)-f\big(u_m(z)\big)\Big]_++\left(\frac{1}{2n}+
\frac{1}{2m}\right)\right)\,dz\\
&\leq T\left(\frac{1}{2n}+\frac{1}{2m}\right)+ \int_0^t\beta\Big(\big\|u_n(z)-u_m(z)\big\|\Big)\,dz.
\end{align*}
Let $M$ come from the definition of slow function (see Subsection \ref{sf}).  Since $\beta$ is a slow function we get
$$\big\|u_n(t)-u_m(t)\big\|\leq T\left(\frac{1}{2n}+\frac{1}{2m}\right)M\to 0\;\; \text{as}\;\; n,m\to\infty.$$
\indent Let $u(t):=\lim_{n\to\infty}u_n(t)$, $t\in [0,T]$. Then $u:[0,T]\to X$ is an integral solution to \eqref{cauchy-accretive-f}. Indeed, for any $(y,v)\in \Gr(A)$, $0\leq s\leq t\leq T$ and  any $n\in\N$,
\begin{equation}\label{pom1}\big\|u_n(t)-y\big\|\leq\big\|u_n(s)-y\big\|+
\int_s^t\big[u_n(z)-y,f_n(u_n(z))-v\big]_+\,dz.\end{equation}
Passing in \eqref{pom1} with $n\to\infty$, using the upper semicontinuity of $[\cdot,\cdot]_+:X\times X\to\R$ and the $\sup$-Fatou lemma we get that
$$\big\|u(t)-y\big\|\leq\big\|u(s)-y\big\|+\int_s^t\big[u(z)-y,f(u(z))-v\big]_+\,dz.$$
This shows that $u$ is an integral solution to \eqref{cauchy-accretive-f} as required. Condition \eqref{slow one sided} immediately implies that $u$ is a unique solution to \eqref{cauchy-accretive-f}.\\
\indent To show the Lipschitz continuity of $u$ starting at $x\in D(x_0,r)\cap D(A)$ let $0\leq t<s\leq T$, $h:=s-t$ and $w(z):=f\big(u(z)\big)$ for $z\in [0,T]$. As before, for $v\in Ax$,
\begin{align*}\big\|u(t+h)-u(t)\big\|&\leq \big\|u(h)-x\big\|+\int_0^t\Big[u(z+h)-u(z),f\big(u(z+h)\big)-f\big(u(z)\big)\Big]_+\,dz\\
&\leq \int_0^h\big\|w(z)-v\big\|\,dz+\int_0^t
\beta\Big(\big\|u(z+h)-u(z)\big\|\Big)\,dz\\
&\leq ch+\int_0^t
\beta\Big(\big\|u(z+h)-u(z)\big\|\Big)\,dz.\end{align*}
Hence $\big\|u(t+h)-u(t)\big\|\leq chM$, where $c:=\sup_{z\in [0,T]}\big\|w(z)-v\big\|$.
\end{proof}

\begin{rem} \label{integral-mild-coincide} {\em (1) Both Proposition \ref{prop lip} and Theorem \ref{main existence} hold true for arbitrary quasi $m$-accretive operators. The first statement follows immediately from Remark \ref{semigroup prop} (1). The second follows from the fact that if $A$ is $\alpha$-$m$-accretive and condition \eqref{slow one sided} holds, then the function $f+\alpha I$ also satisfies \eqref{slow one sided} with the slow function $s\mapsto \beta(s)+\alpha s$ (see Remark \ref{slow rem}.).\\
\indent  (2) Assume that $X$ is reflexive. In view of Proposition \ref{prop lip}, we have that for $x\in D(x_0,r)$ an integral solution $u:[0,T]\to \Omega$ to \eqref{cauchy-accretive-f} is the uniform limit of a sequence of strong solutions to
\eqref{cauchy-accretive-fn}.\\
\indent (3) If above $f:X\to X$ and $\beta(s)=bs$ for some $b>0$, then $-f$ is $b$-accretive (it is even strongly $b$-accretive). Due to the Martin theorem \cite[page 104]{Barbu}, $f$ is in fact $b$-$m$-dissipative. Thus, by \cite[Theorem 3.1]{Barbu}, $A-f$ is quasi $m$-accretive. Therefore, for any $x\in \overline{D(A)}$, $S_{A-f}(\cdot)x:[0,\infty)\to X$ is well-defined. Let $u:[0,\infty)\to X$ be an integral solution to \eqref{cauchy-accretive-f}, i.e. $u=u_A(\cdot;x,w)$, where $w:=f\circ u$. We shall show that $u=S_{A-f}(\cdot)x$. First observe that, due to Remark \ref{semigroup prop} (2)  without loss of generality we may assume that $b=0$.  Assuming that $f$ is strongly dissipative we have that for any $(y,v)\in\Gr(A)$ and $0\leq s<t$
$$\big\|u(t)-y\big\|\leq\big\|u(s)-y\big\|+\int_s^t\big[u(z)-y,f\big(u(z)\big)-v\big]_+\,dz.$$
But \begin{align*}\big[u(z)-y,f\big(u(z)\big)-v\big]_+&\leq \big[u(z)-y,f\big(u(z)\big)-f(y)\big]_++\big[u(z)-y,f(y)-v\big]_+\\ &\leq \big[u(z)-y,f(y)-v\big]_+.\end{align*}
Hence $$\big\|u(t)-y\big\|\leq\big\|u(s)-y\big\|+\int_s^t\Big[u(z)-y,-\big(v-f(y)\big)\Big]_+\,dz.$$ Since $v\in Ay$ was arbitrary, we see that $v-f(y)$ is an {\em arbitrary} element of $(A-f)(y)$. This shows that $u=S_{A-f}(\cdot)x$.
} \end{rem}

We complete the section with the result allowing us to treat an integral solution of \eqref{equat-A+f} as an integral (or even strong) one in a suitable greater phase space (see Remark \ref{how-to-compute} (2)).
\begin{prop}\label{reg-trick}
Suppose that $Y$, $A_Y$ and $V_Y$ are as in Remark {\em \ref{how-to-compute} (1)}. If $u:[0,T]\to X$ is an integral solution of \eqref{equat-A+f} with $u(0)\in D(A)$, then
$u\in W^{1,1}([0,T],Y)$, $u(t)\in D(A_Y)$ and
$$
\dot u(t) \in - A_Y u(t) + f(u(t)),
$$
for a.e. $t\in [0,T]$.
\end{prop}
\begin{proof} Observe that $w:=f\circ u \in L^1([0,T],Y)$ and that,
since the part of $A$ in $X$  is equal to $A$, we have $(I+\lambda A_Y)^{-1} v = (I+\lambda A)^{-1} v \in D(A)$ for any $v\in X$ and $\lambda>0$. By the construction of integral solutions (see \cite[Subsection 4.1]{Barbu}) it is easily seen that $u$ is also an integral solution of
\begin{equation}\label{Y-equ}
    \dot u(t)\in A_Y u(t)+w(t),  \ t\in [0,T].
\end{equation}
Hence, in view of Proposition \ref{prop lip} (c) and the reflexivity of $Y$, we infer that $u\in W^{1,\infty}([0,T],Y)$ and $u$ is a strong solution of \eqref{Y-equ}.
\end{proof}

\section{Differentiation along trajectories and proofs of Theorems \ref{thm-inv1}--\ref{thm-inv2-prime}}\label{sec3}

Let $K\subset X$ be closed and $x\in X$. Given a continuous curve $u:[0,\tau)\to X$, $\tau>0$, such that $u(t)=x$ for some $t\in [0,\tau)$, the Dini derivative $D_+(d_K\circ u)(t)$ measures the rate of changes of the distance of $u(s)$ from $K$ for $s$ in a neighborhood of $t$. It is easy to see for instance that $D_+(d_K\circ u)(t)=0$ if and only if there exist sequences $h_n\to 0^+$ and $v_n\to 0$ such that $d_K\big(u(t+h_n)+h_nv_n\big)\leq d_K(x)$ for all $n\geq 1$, i.e. $u$ is tangential to the set $K^\alpha:=\big\{y\in X\mid d_K(y)\leq\alpha\big\}$, where $\alpha:=d_K(x)$. If $u$ is (right) differentiable at $t$, i.e. $u'_+(t):=\lim_{h\to 0^+}h^{-1}\big(u(t+h)-u(t)\big)$ exists, then $u$ is tangential to $K^\alpha$  if and only if  the vector  $v=u'_+(t)$ is tangent to $K^\alpha$ at $x$, i.e.
$$D_+d_K(x;v):=\liminf_{h\to 0^+}\frac{d_K(x+hv)-\alpha}{h}=0;$$
in other words if and only if $v\in T_{K^\alpha}(x)$.\\
\indent Let $A:D(A)\to X$ be a quasi $m$-accretive operator, let $V:X\to\R$ be a locally Lipschitz function representing $K$ given by \eqref{constr}, let $x\in\overline{D(A)}$ and $v\in X$. Suppose that $u:=u(\cdot;x,v)=S_{A_v}(\cdot)x:[0,\infty)\to X$, where $A_v:=A(\cdot)-v$, is the integral solution to \eqref{cauchy-accretive-w} with $w(\cdot)\equiv v\in X$ (see also Remark \ref{integral-mild-coincide} (1) with regard to the last equality). By the {\em $A$-derivative of $V$ at $x$ in the direction $v$} we mean the Dini type derivative
\begin{equation}\label{A-derivative}
D_AV(x;v):=\liminf_{h\to 0^+}\frac{(V\circ u)(h)-V(x)}{h}=D_+\big(V\circ u_A (\cdot;x,v)\big)(0).\end{equation}
Note that if $x\in D(A)$, then the derivative $D_AV(x;v)$ is finite since, by Remark \ref{dod info} (a), the function $V\circ u$ is Lipschitz around 0. As above $D_AV(x;v)$ measures the rate of growth of $V$ along the integral curve $u=u_A (\cdot;x,v)$. In particular if $D_A V(x;v)> \alpha$, then  there is $\eta>0$ such that $V\big(u(t)\big)> \alpha t+V(x)$ for $0\leq t<\eta$. It is clear again that  if $\dot u_+(0)$ exists, then
\begin{equation}\label{prawa poch}D_AV(x;v)=D_+V\big(x;\dot u_+(0)\big).\end{equation}
Indeed $u(h)=x+hu'_+(0)+o(h)$ when $h\to 0$; hence
$$h^{-1}\Big|\big(V\big(u(h)\big)-V(x)\big)-\big(V\big(x+h\dot u_+(0)\big)-V(x)\big)\Big|\leq h^{-1}\ell\big|o(h)\big|\to 0\;\; \text{as}\;\; h\to 0,$$
where $\ell$ is the Lipschitz constant of $V$ at $x$.

For a general $A$, $x\in\overline{D(A)}$ and $v\in X$
the $A$-derivative $D_AV(x;f(x))$ is not easy to compute. The situation changes under additional assumptions on $X$ if $x\in D(A)$.

\begin{prop}\label{A-derivative-char}
\indent {\em (i)} If $X$ and $X^*$ are uniformly convex, then for any $x\in D(A)$ and $v\in X$,
\begin{equation}\label{uc}D_AV(x;v)=D_+V(x;v-y),\end{equation}
where $y\in Ax$ is such that $y-v=(Ax-v)^0$ is the element of the set $Ax-v$ having minimal norm. The same holds true in an arbitrary Banach space provided $A$ is a linear operator and then $D_AV(x;v)=D_+V(x;v-Ax)$. \\
\indent {\em (ii)} Assume that $A$ is single-valued and both $X$, $X^*$ are uniformly convex or $X$ is an arbitrary Banach space but $A$ is linear. If for $x\in \overline{D(A)}$
\begin{equation}\label{V-inv}
V\big(S_A (t)x\big)\leq V(x)\;\; \text{for all}\;\; t\geq 0,
\end{equation}
then
\begin{equation}\label{ineq-sep}
D_AV(x;v) \leq V^\circ (x;v)\;\; \text{for any}\;\; x\in D(A),\; v\in X,
\end{equation}
where $V^\circ (x;v)$ is the generalized Clarke derivative at $x$ in the direction of $v$ (see \cite[Chapter 2.1]{clarke}).
\end{prop}
\begin{proof}
\indent (i) If $X$ and $X^*$ are uniformly convex, then in view of Section \ref{cp} (b), $u=u_A (\cdot;x,v)$ is a strong solution and $\dot u_+(0)=(-Ax+v)^0$ exists. Hence $D_A V(x;v)=
D_+V\big(x,(-Ax+v)^0\big)$ in view of \eqref{prawa poch}.\\
\indent If $X$ is arbitrary but $A$ is linear, then in view of Remark \ref{dod info} (c) and \eqref{duhamel}
$$
u_A(t;x,v)= e^{-tA}x+\int_0^te^{-(t-s)A}v\,ds,\; t\geq 0.
$$
Hence
$$\dot u_+(0)=\lim_{t\to 0^+}\frac{e^{-tA}x-x}{t}+\frac{1}{t}\int_0^te^{-sA}v\,ds=-Ax+v.$$
and, consequently, $D_AV(x;v)=D_+V(x;v-Ax)$.\\
\indent (ii) By (i),
$$D_AV(x;v)=D_+V(x;v-Ax)\leq V^\circ(x;v)+D_+V(x;-Ax),$$
where $y\in Ax$ is as in \eqref{uc}.
It is enough to see that $D_A V(x;0)\leq 0$ in view of \eqref{V-inv}. To see this observe that if $u:=S_A(\cdot)x$, then $\dot u_+ (0) = \lim_{t\to 0^+}h^{-1}\big(S_A(t)x-x\big)=-Ax$ exists.
Hence by \eqref{V-inv}
$$D_+V(x,-Ax)=D_AV(x,0)=\liminf_{h\to 0^+}\frac{V\big(S_A(h)x\big)-V(x)}{h}\leq 0.\eqno\qedhere$$
\end{proof}

\indent The following property enables to study the behavior of $V$ along an integral solution curve to \eqref{equat-A+f} without a prior knowledge of this solution.

\begin{prop}\label{diff-superposition} Let $x_0\in \overline{D(A)}$ and $u:[0,\tau)\to X$ be an integral solution to \eqref{equat-A+f}.  Then
$$D_+(V\circ u)(t)= D_AV\Big(u(t);f\big(u(t)\big)\Big),\;\; t\in [0,\tau).$$
\end{prop}
\begin{proof} Fix $t\in [0,\tau)$, let $x:=u(t)$ and $v:=f\big(u(t)\big)$. Recall $D_AV(x;v) =D_+\big(V\circ u_A (\cdot;x,v)\big)(0)$.
In view of Remark \ref{semigroup prop} (3) $u(t+h)=u(h;x,w)$, where $w(h):=f\big(u(t+h)\big)$, for $h\in [0,T-t)$. By \eqref{semi-equality} we have
\begin{gather*}\big|V\big(u(t+h)\big)-V\big(u_A (h;x,v)\big)\big|\leq\ell\big\|u(t+h)-u_A (h;x,v)\big\|=\\
\ell\big\|u_A (h;x,w)-u_A (h;x,v)\big\|\leq\ell e^{\alpha h}\int_0^h\big\|w(s)-v\big\|\,ds,\end{gather*}
where $\ell$ is the Lipschitz constant of $V$ around $x$. Since
$$h^{-1}\ell e^{\alpha h}\int_0^h\big\|w(s)-v\big\|\,ds\to 0\; \text{as}\; t\to 0^+$$ and
$$V\big(u(t+h)\big)-V\big(u(t)\big)=V\big(u(t+h)\big)-V\big(u_A (h;x,v)\big)+V\big(u_A (h;x,v)\big)-V(x)$$
this yields that $D_+(V\circ u)(t)=D_AV(x;v)=D_AV\Big(u(t);f\big(u(t)\big)\Big)$ as required. \end{proof}

This gives us immediately the invariance criteria mentioned in Introduction.

\noindent {\em Proof of Theorem \ref{thm-inv1}.} Suppose to the contrary that there is an integral solution $u:[0,\tau)\to X$ that leaves $K$, i.e. there is $T\in (0,\tau)$ such that $V\big(u(T)\big)>0$. Let
$$\bar t:=\sup\big\{t\in [0,T]\mid V\big(u(t)\big)\leq 0\big\}.$$
Clearly $\bar t<T$, $V\big(u(\bar t)\big)=0$ and
$$V\big(u(t)\big)>0\;\; \text{for all}\;\; t\in (\bar t,T].$$
In order to simplify the notation and without loss of generality we may assume that $\bar t=0$ and $u(t)\in U\big(u(0)\big)$ for $t\in [0,T]$.
In view of Proposition \ref{diff-superposition}, for any $t\in (0,T]$
\begin{equation}\label{dense-use}
D_+(V\circ u)(t)=D_AV\Big(u(t);f\big(u(t)\big)\Big)\leq \omega\big(V\big(u(t)\big)\big).
\end{equation}
This however, by Lemma \ref{Perron}, implies that $V\big(u(t)\big)=0$ for all $t\in [0,T]$ since $\omega$ is a uniqueness function. \hfill $\square$

\noindent {\em Proof of Theorem \ref{thm-inv1-prime}.}
It is sufficient to slightly modify the proof of Theorem \ref{thm-inv1}.
Namely, use condition \eqref{dense-prime} from the assumption
to have the counterpart of the relation \eqref{dense-use} directly, i.e.
$$
D_{+}(V\circ u) (t) \leq \omega \big(V\big(u(t\big)\big) \ \mbox{ for a.e. } t\in (0,T].
$$
This allows us to use Lemma \ref{Perron} to complete the proof. \hfill $\square$

\noindent {\em Proof of Theorem \ref{thm-inv2}.} Suppose to the contrary that there is an integral solution $u_0:[0,\tau)\to X$ of \eqref{equat-A+f} that leaves $K$, i.e. there is $T\in (0,\tau)$ such that $V(u_0(T))>0$.
As in above proofs, we may assume that $x_0=u_0(0)\in \part K$ and $u_0(t)\in U:=U\big(u_0(0)\big)$ for all $t\in [0,T]$.

\indent Without loss of generality we may also assume that  $V$ is Lipschitz on $D(x_0,\delta)$, where $\delta>0$ is given by assumption, with the constant $\ell>0$. Moreover, take $0<r<\delta$ and a sequence $(x_n)$ such that $x_n\to x_0$, $x_n\in D(A)\cap (D(x_0,r)\setminus K_V)$. In view of Theorem \ref{main existence}, there is a time $T'\leq T$ such that for all (sufficiently large if necessary) $n\geq 1$ there is an integral solution $u_n:[0,T']\to\Omega$ to the problem
\begin{align}\label{cauchy-accretive-f-n}
&\left\{ \begin{array}{l}
\dot u  \in -Au+f(u)\\
u(0)=x_n\in D(A),
\end{array} \right.
\end{align}
With no loss of generality we may assume that $T'=T$. It is clear that the sequence $(u_n)$ converges uniformly to $u_0$ on $[0,T]$ since $u_0$ is a unique integral solution to \eqref{equat-A+f}. According to Proposition \ref{prop lip} (b), (c) and
Remark \ref{integral-mild-coincide} (1) for each $n\geq 1$ the function $u_n$ is locally Lipschitz and, due to the reflexivity of $X$, $u_n$ is a strong solution. Therefore there is a full-measure set $S\subset [0,T]$ such that  $u_n(t)\in D(A)$ for all $n\geq N$ and $t\in S$, and $V\circ u_n$ is differentiable at $t$.\\
\indent Take a small $0<\eps<T$. Since $(u_n)$ converges uniformly to $u_0$ we may assume that $u_n(t)\in U\setminus K$ for all $n\geq 1$ and $t\in [\eps,T]$. Let $n\geq 1$ and $z\in S$, $\eps\leq z<T$. Then, in view of Remark \ref{dod info} (d) and \eqref{wzorek},
$u_n(z+\cdot)=u_A (\cdot;u_n(z),w_n)$, where $w_n=f\big(u_n(z+\cdot)\big)$, and
\begin{align*} V\big(u_n(z+h)\big)-V\big(u_n(z)\big)=&\;V\big(u_n(z+h)\big)-V(u_A (h;u_n(z),v_n)
+V(u_A (h;u_n(z),v_n)-V\big(u_n(z)\big),
\end{align*}
where $v_n:=f\big(u_n(z)\big)$. Note that
\begin{align*}\big|V\big(u_n(z+h)\big)-V(u(h;u_n(z),v_n)\big|\leq \ell\big\|u_n(z+h)-u(h;u_n(z),v_n)\big\|\\
=\ell\big\|u(h;u_n(z),w_n)-u(h;u_n(z),v_n)\big\|\leq\ell\int_0^h\big\|w_n(\xi)-v_n\big\|\,d\xi,
\end{align*}
where, as before, $\ell$ is the Lipschitz constant of $V$ around $x_0$.
Therefore
\begin{align*}
(V\circ u_n)'(z)=&\,\lim_{h\to 0^+}\frac{V\big(u_n(z+h)\big)-V\big(u_n(z)\big)}{h}
\\\leq&\,\lim_{h\to 0^+}\frac{\ell}{h}\int_0^h\big\|w_n(\xi)-v_n\big\|\,d\xi+
\liminf_{h\to 0^+}\frac{V(u(h;u_n(z),v_n)-V\big(u_n(z)\big)}{h}\\
=&\,
D_A V\Big(u_n(z);f\big(u_n(z)\big)\Big)\leq \omega\big(V\big( u_n(z)\big)\big).
\end{align*}
This implies that for all $t\in [\eps,T]$
\begin{equation}\label{dense-bis}
V\big(u_n(t)\big)=V\big(u_n(\eps)\big)+\int_\eps^t(V\circ u_n)'(z)\,dz\leq V\big(u_n(\eps)\big)+\int_\eps^t\omega\Big(V\big(u_n(z)\big)\Big)\,dz.
\end{equation}
Passing with $\eps\to 0^+$ and with $n\to\infty$ we get for all $t\in [0,T]$
$$
V\big(u_0(t)\big)\leq \int_0^t\omega\Big(V\big(u_0(z)\big)\Big)\,dz.
$$
Since $\omega$ is nondecreasing, we have
\[\left(\int_0^t\omega\Big(V\big(u_0(z)\big)\Big)\,dz\right)' = \omega\Big(V\big(u_0(t)\big)\Big)\leq \omega\left(\int_0^t\omega\Big(V\big(u_0(z)\big)\Big)\,dz\right),\]
which implies that $\int_0^\cdot\omega\Big(V\big(u_0(z)\big)\Big)\,dz\equiv 0$ and, in consequence, due to Lemma \ref{Perron}, we get $V\big(u_0(t)\big) = 0$ for all $ t\in [0,T] $: a contradiction.\hfill $\square$

\noindent {\em Proof of Theorem \ref{thm-inv2-prime}.}
A slight modification of the proof of Theorem \ref{thm-inv2} is sufficient. Without the reflexivity of $X$ we still know that $u_n$ are locally Lipschitz in view of Proposition \ref{prop lip} (b) and that $V\circ u_n$ is locally Lipschitz and, hence, a.e. differentiable. In particular, for a.e. $t\in [\varepsilon, T]$,
$$
(V\circ u_n)'(t)=D_+ (V\circ u_n)(t).
$$
Applying assumption (iii)' one has \eqref{dense-bis}, which allows us to complete the proof by following the steps of the previous one. \hfill $\square$

\section{Strict invariance. Proofs of theorems \ref{thm-inv-kv0-1} -- \ref{thm-strinv-2}}\label{sec-str-inv}

Recall that a closed $K\subset\Omega\cap\overline{D(A)}$ is strictly invariant with respect to \eqref{equat-A+f} if all solutions starting at $x_0\in K$ stay for $t\in (0,\tau_u)$ in the interior $\wew K$ of $K$.
In particular, a {\em necessary condition for the strict invariance} is that  $\wew K$ is invariant because solutions cannot return from $\wew K$ to $\part K$. At the beginning of the discussion we state conditions implying this necessary condition, i.e. the invariance of $\wew K$. Recall that $K=K_V$ for some $V$ -- see \eqref{constr}, and assume that $K_V^0=\wew K$.

\noindent {\em Proof of Theorem \ref{thm-inv-kv0-1}.}
Suppose to the contrary that there is an integral solution $u:[0,\tau)\to X$ of \eqref{equat-A+f} starting at $x\in K_V^0$ and leaving it, i.e. $V\big(u(t)\big)=0$ for some $t\in (0,\tau)$. Let
$$T:=\inf\big\{t\in [0,\tau)\mid V\big(u(t)\big)= 0\big\}.$$
Clearly $T>0$, $V\big(u(t)\big)<0\;\; \text{for all}\;\; t\in [0,T)$, and $u(T)\in \partial K$. There is  $\bar t\in [0,T)$ such that $u(t)\in U\big(u(T)\big)$ for $t\in [\bar t,T]$ and without loss of generality we may suppose that $\bar t=0$.
In view of Proposition \ref{diff-superposition} and  assumption (iv), for any $t\in [0,T)$,
\begin{equation}\label{nier}D_+(V\circ u)(t)=D_AV\Big(u(t);f\big(u(t)\big)\Big)\leq \omega\big(-V\big( u(t)\big)\big).\end{equation}
By Lemma \ref{Perron} (2) applied to $v:[-T,0]\to X$ given by $v(s):=(V\circ u)(s+T)$ for $-T\leq s\leq 0$, we see that \eqref{nier} implies that $V\big(u(t)\big)=0$ for all $t\in [0,T]$; a contradiction.\hfill $\square$

\noindent {\em Proof of Theorem \ref{thm-inv-kv0-2}.} Suppose to the contrary that there is an integral solution $u_0:[0,\tau)\to X$ of \eqref{equat-A+f} and $0<T<\tau$ such that   $V\big(u_0(T)\big)=0$ and $V\big(u_0(t)\big)<0\;\; \text{for all}\;\; t\in [0,T)$.
Without loss of generality we may assume that $u_0(t)\in U:=U\big(z_0\big)$ for all $t\in [0,T]$, where $z_0=u_0(T)$. Clearly $z_0\in\part K$.\\
\indent There is $\delta>0$ such that  $V$ is Lipschitz on $D(z_0,\delta) \subset U$ with the constant $\ell>0$. Take $0<r<\delta$. In view of Theorem \ref{main existence} there  is $0<T'<\tau-T$ such that all solutions to \eqref{equat-A+f} starting in $D(z_0,r)$ are defined on $[0,T']$. Take $t_0\in [0,T)$ such that $u_0(t)\in D(z_0,r)$ for every $t\in (t_0,T]$ and $T'':=T-t_0\leq T'$. Now, take a sequence $(x_n)$ such that $x_n\to x_0:=u_0(t_0)$, $x_n\in D(A)\cap \big(D(z_0,r)\cap K_V^0\big)$. For all (sufficiently large if necessary) $n\geq 1$ there is an integral solution $u_n:[0,T'']\to\Omega$ to problem \eqref{cauchy-accretive-f-n}.\\
\indent It is clear that the sequence $(u_n)$ converges uniformly to $\bar u:=u_0(t_0 +\cdot)$ on $[0,T'']$ since $u_0$ is a unique integral solution to \eqref{equat-A+f}. Each of solutions $u_n$ is Lipschitz continuous and $u_n(t)\in D(A)$ for a.a. $t\in [0,T'']$. Hence there is a set $S\subset [0,T]$ of full measure such that for all $n\geq N$ $u_n(t)\in D(A)$ and $V\circ u_n$ is differentiable at $t\in S$.\\
\indent Take a small $0<\eps<T''$. Since $(u_n)$ converges uniformly to $\bar u$, we may assume that $u_n(t)\in U\cap K_V^0$ for all $n\geq 1$ and $t\in [0,T''-\eps]$. Let $n\geq 1$ and $z\in S$, $0< z\leq T''-\eps$. Then, in view of Remark \ref{dod info} and  \eqref{wzorek},
$u_n(z+\cdot)=u_A \big(\cdot;u_n(z),w_n\big)$, where $w_n=f\big(u_n(z+\cdot)\big)$, and we can repeat arguments from the proof of Theorem \ref{thm-inv2} to obtain
\[(V\circ u_n)'(z)\leq
D_AV\Big(u_n(z);f\big(u_n(z)\big)\Big)\leq \omega\big(-V\big(u_n(z)\big)\big).
\]
This implies that for all $t\in [0,T''-\eps]$
$$V\big(u_n(t)\big)=V\big(u_n(T''-\eps)\big)-\int_t^{T''-\eps}(V\circ u_n)'(z)\,dz\geq V\big(u_n(T''-\eps)\big)-\int_t^{T''-\eps}\omega\Big(-V\big(u_n(z)\big)\Big)\,dz.$$
Passing with $\eps$ to $0$ and with $n$ to infinity we get for all $t\in [0,T'']$
$$-V\big(\bar u(t)\big)\leq \int_0^t\omega\Big(-V\big(\bar u(z)\big)\Big)\,dz$$
i.e. $0\leq -V\big(u_0(t)\big)\leq \int_0^t\omega\Big(-V\big(u_0(z)\big)\Big)\,dz$ for all $t\in [t_0,T]$.
As in the proof of the previous theorem this implies that $V(u_0(t))=0$ for all $t\in [t_0,T]$: a contradiction.\hfill $\square$

Finally we are in a position to prove the strict invariance results.

\noindent {\em Proofs of Theorem \ref{thm-strinv-1} and Theorem \ref{thm-strinv-2}.} From assumptions \eqref{dense2} and \eqref{strict-inward} (resp. \eqref{weak dense2} and \eqref{strict-inward}) and conclusions of theorems \ref{thm-inv1} and \ref{thm-inv-kv0-1} (resp. \ref{thm-inv2} and \ref{thm-inv-kv0-2}) it follows that both $K=K_V$ and $K_V^0$ are invariant. Suppose that there is an integral solution $u$ of \eqref{equat-A+f} such that $u(t)\in \part K$ for some interval $[0,T]$, $T>0$. Then $(V\circ u)(t)=0$ for all $t\in [0,T]$, which implies (see Proposition \ref{diff-superposition}) that $$D_AV\Big(u(t);f\big(u(t)\big)\Big)=D_+(V\circ u)(t)=0 ,\;\; \mbox{for}\;\; t\in [0,T);$$
a contradiction.\hfill $\square$

\section{Applications}\label{applications}

\subsection{Impulsive differential equations with state-de\-pen\-dent impulses}
Problems with state-de\-pen\-dent impulses of the form
\begin{equation}\label{problem}
\left\{ \begin{array}{ll}
    \dot{y}(t)\in F\big(t,y(t)\big),\;\; t\in[0,T],\;t\neq \tau_j\big(y(t)\big),\;j=1,\ldots,k, \\
    y(0)=x_0,  \\
    y(t^+)=y(t)+I_j\big(y(t)\big)\;\; \text{for}\;\; t=\tau_j\big(y(t)\big),\; j=1,\ldots,k,
\end{array}\right.
\end{equation}
where $T>0$, $F:[0,T]\times X\multimap X$ is a set-valued dynamics, for $j=1,...,k$, $\tau_j:X\to (0,T)$ is a {\em barrier function} and $I_j:X\to X$  an {\em impulse function}, meet a considerable interest recently. In order to  characterize the suitable function space where solutions can be considered, one looks for sufficient conditions implying that every trajectory of \eqref{problem} meets a barrier $\Gamma_j=\Gr(\tau_j)$ {\em exactly} once. Note that if the global existence is achieved, then each barrier is hit {\em at least} once. One however demands that after the $j$-th jump a solution stays in the epigraph $\Epi(\tau_j)$ of $\tau_j$, i.e. it immediately enters its interior and does not return to $\Gamma_j$; in other words one needs conditions implying that epigraph $\Epi(\tau_j)$, $j=1,...,k$, is strictly invariant. Results from Section \ref{sec-str-inv} do fit well to this problem if $F(t,y)=-Ay+f(t,y)$, where $A:D(A)\multimap X$ is an $m$-accretive operator.\\
\indent Let us consider the following  problem
\begin{equation}\label{impuls}\begin{cases}
\dot u\in-Au+f(t,u),\;\; t\in [0,T],\\
u(0)=x\in \overline{D(A)},\end{cases}
\end{equation}
where $f:\R\times X\to X$ is continuous. Let $\tau:X\to\R$ be a locally Lipschitz barrier function. By a {\em solution} to \eqref{impuls} we understand an integral solution $u:[0,T]\to X$, $T>0$, to  \eqref{cauchy-accretive-w} with $w=f\big(\cdot,u(\cdot)\big)$.

\begin{thm}\label{thm-strinv-ext}
Assume that for every $(z,\theta)\in \Gr(\tau)$ there are a neighborhood $U=U(z,\theta)$ of $(z,\theta)$ and a uniqueness function $\omega$ such that
\begin{gather}\label{dense3}D_A\tau\big(x;f(t,x)\big)\leq\omega\Big(\big|\tau(x)-t\big|\Big)+1\;\; \text{for}\;\; (t,x)\in U,\\
\label{strict-inward1}
D_A\tau\big(x;f(t,x)\big)<1\;\; \text{for}\;\; (t,x)\in U\cap \Gr(\tau).\end{gather}
If $u:[0,T]\to X$ is a solution to \eqref{impuls} and $\tau(x_0)\leq t$, then $\tau\big(u(h)\big)<t+h$ for any $0<h\leq T$, i.e. $\big(t+h,u(h)\big)\in\Epi(\tau)$ for $h\in (0,T]$.
\end{thm}

\begin{proof} Define $\A:D(\A)\to \X:=\R\times X$ by $\A(t,u):=(0,Au)$ for $(t,u)\in D(\A):=\R\times D(A)$ and $\F:\X\to\X$ by $\F(t,x):=\big(1,f(t,x)\big)$ for $(t,x)\in\X$; $\X$ is a Banach space with the norm $\big\|(t,x)\big\|:=|t|+\|x\|$ for $(t,x)\in\X$.
It is immediate to see that $\A$ is $m$-accretive and $\F$ is continuous. A simple calculation shows that a continuous function $\u:[0,T]\to\X$ is an integral solution to the problem
\begin{equation}\label{impuls1}\begin{cases}\dot\u\in-\A\u+\w,\\
\u(0)=(t,x),\end{cases}\end{equation}
with $(t,x)\in\X$ and $\w=(1,w)$, where $w\in L^1\big([0,T],X\big)$, if and only if $\u(h)=\big(t+h,u(h)\big)$ for $h\in [0,T]$, where $u:[0,T]\to X$ is an integral  solution to \eqref{cauchy-accretive-w}. In particular, $u:[0,T]\to X$ is a solution of \eqref{impuls} if and only if $\u(h)=\big(t+h,u(h)\big)$ is a solution to \eqref{impuls1} with $\w=\F\circ\u$. \\
\indent Let $V:\X\to\R$ be given by $V(t,x):=\tau(x)-t$, $(t,x)\in\X$.
Clearly $V$ is locally Lipschitz, $\K:=\big\{(t,x)\in\X\mid V(t,x)\leq 0\big\}=\Epi(\tau)$, $\wew\K=\big\{(t,x)\in\X\mid V(t,x)<0\big\}$ and $\part\K=\Gr(\tau)$.\\
\indent For a fixed $(t,x)\in\X$  let $\bar\u$ denote the solution to \eqref{impuls1} with $\w(\cdot)\equiv \F(t,x)$. Then $\bar\u=\big(t+\cdot,\tilde u(\cdot)\big)$ on $[0,T]$, where $\tilde u=u_A\big(\cdot;x,f(t,x)\big)$. Therefore
\begin{align*}
D_\A V\big((t,x);\F(t,x)\big)&=\liminf_{h\to 0^+}\frac{V\big(\bar\u(h)\big)-V(t,x)}{h}=\liminf_{h\to 0^+}\frac{\tau\big(u(h)\big)-t-h-\tau(x)+t}{h}\\
&=\liminf_{h\to 0^+}\frac{\tau\big(u(h)\big)-\tau(x)}{h}-1=D_A\tau\big(x;f(t,x)\big)-1.
\end{align*}
Hence, by \eqref{dense3} and \eqref{strict-inward1},
\begin{gather*}D_\A V\big((t,x);\F(t,x)\big)\leq \omega\big(|V(t,x)|\big)\;\; \text{for each}\;\; (t,x)\in U,\\
D_\A V\big((t,x);\F(t,x)\big)<0\;\; \text{for each}\;\; (t,x)\in U\cap \Gr(\tau).\end{gather*}
In view of Theorem \ref{thm-strinv-1}, the set $\K$ is strictly invariant with respect to \eqref{impuls1}. This completes the proof.
\end{proof}

\begin{rem}
\begin{em}
Assumptions \eqref{dense3} and \eqref{strict-inward1} allow to consider nonsmooth barriers and deal with integral (not strong) solutions in contrast to many other papers, see e.g. \cite{hptt}. Note also that in our approach the operator $A$ may be nonlinear and solutions need not be mild as, e.g., in \cite{bcgr}.
\end{em}

\end{rem}

\begin{cor}
Assume that $A$ and $f$ satisfy the assumptions of Theorem \ref{thm-strinv-ext}, and let $\tau_j:H\to (0,\infty)$ be locally Lipschitz functions such that \eqref{dense3} and \eqref{strict-inward1} hold for $\tau_j$ instead of $\tau$, for every $j=1,\ldots,k$. We also assume standard conditions on barriers:
\begin{enumerate}[\em ($\tau$1)]
\item $0<\tau_j(x)<\tau_{j+1}(x)$ for each $x\in X$ and $j=1,\ldots,k$,
\item $\tau_j\big(x+I_j(x)\big)\leq\tau_j(x)<\tau_{j+1}\big(x+I_j(x)\big)$ for each $x\in X$ and $j=1,\ldots,k$.
\end{enumerate}
Then any solution to \eqref{impuls} meets each barrier $\Gamma_j$ at most once.\hfill $\square$
\end{cor}

\subsection{Obstacle problem for equations with one dimensional $p$-Laplace operator}
Consider the following nonlinear problem
\begin{equation}\label{p-laplace-1-dim}
\left\{
\begin{array}{ll}
u_t =  \Delta_p u + f(x,u),\;&x\in (0,l), \, t\in [0,T],\\
u(0,t) = u(l,t) = 0,\;&t\in [0,T],\\
\end{array}
\right.
\end{equation}
where $l>0$, $T>0$, $\Delta_p u:=\big(|u_x|^{p-2}u_x\big)_x$, $p\geq 2$, is the so-called $p$-Laplacian and  $f:[0,l]\times \mathbb{R}\to \mathbb{R}$. Suppose that functions $m, M:[0,l]\to\R$ such that
\begin{equation}\label{dot m}m\leq M,\;\; m(0)\leq 0\leq M(0)\;\; \text{and}\;\; m(l)\leq 0\leq M(l)\end{equation} represent the obstacles. In the so-called {\em obstacle problem} we look for conditions on $f$, $m$ and $M$ implying that for any continuous $u_0:[0,l]\to\R$ such that $u_0 (0)= u_0 (l)=0$ and $m(x)\leq u_0(x)\leq M(x)$ for  $x\in [0,l]$, all solutions of \eqref{p-laplace-1-dim} starting at $u_0$ satisfy
\begin{equation}\label{obstacle-problem-p-lap}
m(x)\leq u(x,t)\leq M(x)\;\; \text{for all}\;\; x\in [0,l],\; t\in [0,T].\end{equation}
\indent Let $X=C_0[0,l]=\big\{ u \in C[0,l]\mid u(0)=u(l)=0\big\}$, where $C[0,l]$ stands for the space of continuous functions $u:[0,l]\to\R$. Clearly $X$ endowed with the sup-norm $\|\cdot\|_\infty$ is a Banach space. Let $A:D(A)\to X$ be given by $Au:=-\big(|u'|^{p-2}u'\big)'$ for $u\in D(A)$, where
$$D(A):=\Big\{u\in X\cap C^1(0,l)\, \mid \big(|u'|^{p-2}u'\big)' \in X\Big\}.
$$
By $C^k(0,l), k\geq 1$, we mean the space of functions from $C[0,l]$ with continuous $k$-th derivative on $(0,l)$.
$D(A)$ is dense and the operator $A$ is $m$-accretive (see e.g. \cite[Lemma 6.1]{Cw-Kry-2009}). In addition to the above assumptions suppose that:
\begin{enumerate}
\item $f:[0,l]\times\R\to\R$ is continuous, $f(\cdot,0)\equiv 0$ and $f(x,\cdot)$ is locally Lipschitz continuous uniformly with respect to $x\in [0,l]$, i.e. for any $s\in\mathbb{R}$ there are $L>0$ and $\delta>0$ such that
\begin{equation}\label{lips}\big|f(x,s_1)-f(x,s_2)\big|\leq L|s_1-s_2|\;\; \text{for all}\;\; x\in [0,l]\;\;\text{and}\;\; s_1,s_2\in (s-\delta,s+\delta);\end{equation}
\item $m, M\in C[0,l]\cap C^2(0,l)$, $m$ is a subsolution and $M$ is a supersolution
of the stationary problem related to \eqref{p-laplace-1-dim}, i.e.
\begin{equation}\label{low-p-lap-ob-cndtn}
- \Delta_p m(x) \leq f\big(x,m(x)\big), \ x\in (0,l)\;\; \text{ and } -\Delta_pM(x)\geq f\big(x,M(x)\big)\, ,\;\;x\in (0,l).
\end{equation}
\end{enumerate}
Condition (1) implies that the Nemytskii operator $F:X\to X$ given by $F(u)(x):= f\big(x,u(x)\big)$, $u\in X$ and $x\in [0,l]$, is well-defined and locally Lipschitz.\\
\indent By a {\em solution} to \eqref{p-laplace-1-dim} on  $[0,T]$ we understand an integral solution  $u:[0,T]\to X$ of the problem
\begin{equation}\label{abs-p-lap-eq}
\dot u = - Au+ F(u), \ t\in [0,T].
\end{equation}
It is clear that a solution $u$ satisfies condition \eqref{obstacle-problem-p-lap} if and only if $u(t)\in K_m\cap K_M$ for all $t\in [0,T]$, where
$$ K_m:=\big\{u \in X\mid u\geq m\; \text{on}\; [0,l] \big\},\;\; K_M:=\big\{ u \in X\mid u\leq M\; \text{on}\; [0,l] \big\}.$$
Clearly
$$K_m=\big\{u\in X\mid V_m(u)\leq 0\big\},\;\; K_M=\big\{u\in X\mid V_M(u)\leq 0\big\},$$
where  $V_m, V_M: X\to \mathbb{R}$  are given by
\begin{equation}\label{V-m-M-formulas}
V_m(u):= \frac{1}{2}\int_{0}^{l} (u-m)_-^2\, d x,\;\; V_M (u):= \frac{1}{2} \int_{0}^{l} (u-M)_+^2\, dx\;\; \text{for}\;\; u\in X.
\end{equation}
\begin{rem} {\em It can be easily verified that $d_{K_m}(u)=\big\|(u-m)_-\|_\infty$ and $d_{K_M} (u) =\big\|(u-M)_+ \big\|_{\infty}$. Observe that neither $V_m=d_{K_m}$ nor $V_M=d_{K_M}$.}
\end{rem}
In order to verify condition \eqref{dense-prime} of Theorem \ref{thm-inv1-prime} we shall follow the idea from Remark \ref{how-to-compute} (2) with $Y=L^2(0,l)$ and the $L^2$-realization of the $p$-Laplace operator. It is known (see e.g. \cite[Th. 3.5]{Cw-Mac}) that, for any $u_0\in X$, any integral solution of \eqref{abs-p-lap-eq} has the following properties
$$
u\in C\big([0,T],X\big) \cap C\big((0,T], W_{0}^{1,p}(0,l)\big) \cap W_{loc}^{1,2}\big((0,T],L^2[0,l]\big),
$$
$\Delta_p u(t)\in L^2[0,l]$ and
\begin{equation}\label{p-lap-reg-sol}
\dot u (t) = \Delta_p u(t) + F\big(u(t)\big)\;\; \text{for a.e.}\;\; t\in [0,T],
\end{equation}
i.e.
$$
\dot u(t) = - A_{L^2} u(t) + F(u(t)), \ \mbox{ for a.e. } \ t\in [0,T],
$$
where $A_{L^2}:D(A_{L^2})\to L^2[0,l]$ is given by
$$
A_{L^2} u:= -\Delta_p u, \  \ D(A_{L^2}):=\{ u\in W^{1,p}(0,l)\mid \Delta_p u \in L^2[0,l] \}.
$$
It is well known that $A_{L^2}$ is $m$-accretive (see e.g. \cite[Lem. 3.4]{Cw-Mac}) and clearly the part of $A_{L^2}$ in $X$ is equal to $A$.
\begin{lm}\label{p-sol-der}
If $u:[0,T]\to X$ is a solution to \eqref{p-laplace-1-dim}, then $V_m\circ u$ \ and  $V_M\circ u$ \ are absolutely continuous on compact subsets of $(0,T]$ and for a.e. $t\in [0,T]$, $\Delta_p u(t) \in L^2[0,l]$  and
\begin{align*}
(V_m\circ u)'(t)&= -\int_{0}^{l} \big(u(t)-m\big)_-\big(\Delta_p u(t) +f(x,u(t))\big)\, dx,\\
(V_M\circ u)'(t)&= \int_{0}^{l} (u(t)-M)_+\big(\Delta_p u(t) +f(x,u(t))\big)\, dx.
\end{align*}
\end{lm}
\begin{proof}
Define $\widetilde V_m, \widetilde V_M:L^2[0,l]\to \mathbb
{R}$ by \eqref{V-m-M-formulas}. It can be easily proved that $\widetilde V_m$ and $\widetilde V_M$ are continuously differentiable and that, for any $w,v\in L^2[0,l]$,
\begin{equation}\label{d-1-dim-der}
\big\la\widetilde V_m'(w),v\big\ra=-\int_{0}^{l} (w-m)_-v\,dx,\;\; \big\la\widetilde V_M'(u),v\big\ra=\int_{0}^{l} (w-M)_+ v\,dx;
\end{equation}
recall that $\la\cdot,\cdot\ra$ stands here for the conjugation duality in $L^2[0,l]$. These formulae follow from \cite[Section 20.6]{Krasno} and \cite[Section 7.4]{Gil}.
Since $X$ is continuously embedded in $L^2[0,l]$, $V_m\big(u(t)\big)=\widetilde V_m\big(u(t))\big)$ and $V_M\big(u(t)\big)=\widetilde V_M\big(u(t)\big)$ for all $t\in [0,T]$, we infer by the chain rule that, for a.e. $t\in [0,T]$,
$$
(V_m\circ u)'(t)=\big\la\widetilde V_m'\big(u(t)\big),\dot u(t)\big\ra\;\; \text{and}\;\;
 (V_M\circ u)'(t)=\big\la\widetilde V_M'\big(u(t)\big),\dot u(t)\big\ra.
$$
This along with \eqref{d-1-dim-der} and \eqref{p-lap-reg-sol} yields the assertion.
\end{proof}
\begin{prop}\label{p-lap-dir-der}
If conditions \eqref{low-p-lap-ob-cndtn} hold, then there exist $C>0$ and $\delta>0$ such that
\begin{align}
\label{fir}-&\int_{0}^{l} (u-m)_-\big(\Delta_p u + f(x,u)\big) \, dx \leq C\,V_m(u)\;\; \text{for any}\;\; u\in D(A_{L^2})\cap \big(B(K_m, \delta)\setminus K_m\big),\\
\label{sec}&\int_{0}^{l} (u-M)_+\big(\Delta_p u+f(x,u)\big)\, dx \leq C\,V_M(u)\;\; \text{for any}\;\; u\in D(A_{L^2}) \cap \big(B(K_M,\delta)\setminus K_M\big).
\end{align}
\end{prop}
\begin{proof} In view of \eqref{lips} and \eqref{low-p-lap-ob-cndtn}, there are $L>0$ and $\delta>0$ such that, for any  $x\in (0,l)$ and $s\in [-\delta,0)$,
\begin{equation}
\begin{split}
&s\Big(\Delta_p m(x)+f\big(x,m(x)+s\big)\Big) \leq s\Big(\Delta_p m(x)+f\big(x,m(x)\big)\Big)\\+&s\Big(f\big(x,m(x)+s\big) - f\big(x,m(x)\big)\Big)
\leq |s|\Big|f\big(x,m(x)+s\big)-f\big(x,m(x)\big)\Big| \leq Ls^2.\end{split}  \label{sub-sol-lip}
\end{equation}
Take $u\in D(A_{L^2})\cap\big(B(K_m,\delta)\setminus K_m\big)$. We then have $0<m(x)-u(x) < \delta$  for all $x\in\big\{x\in [0,l]\mid u(x)<m(x)\big\}\neq\emptyset$.
Observe that
\begin{align*}
&-\int_{0}^{l}(u-m)_-\big(\Delta_p u + f(x,u)\big) \, dx  = \\
& =-\int_{0}^{l}  (u-m)_-(\Delta_p u- \Delta_p m)\, dx   - \int_{0}^{l}(u-m)_-\big(\Delta_p m+ f(x,u)\big) \, dx.
\end{align*}
By \cite[Lemme 5.1]{Glo} for $s_1, s_2\in\R$ one has $\big(|s_1|^{p-2}s_1-|s_2|^{p-2}s_2\big)\cdot (s_1-s_2)\geq |s_1-s_2|^p$. In view of \eqref{dot m}, $(u-m)_- \in X$, so by integrating by parts, we get
\begin{align*}
&\int_{0}^{l}(u-m)_-(\Delta_p u- \Delta_p m)\, dx =-\int_{0}^{l}\big((u-m)_-\big)'\big(|u'|^{p-2} u' - |m'|^{p-2}m'\big) \, dx\\
&= \int_{\{u<m \}}(u'-m')\big(|u'|^{p-2} u' - |m'|^{p-2}m'\big)\,dx\geq 0.
\end{align*}
By applying \eqref{sub-sol-lip} we get
\begin{align*}
& -\int_{0}^{l}(u-m)_-\big(\Delta_p u + f(x,u)\big) \, dx  \leq - \int_{0}^{l}(u-m)_-\big( \Delta_p m+ f(x,u)\big) \, dx\\
= &\int_{\{u<m\}} (u-m)\big(\Delta_p\, m+ f(x,u)\big) \, dx \leq  L\int_{\{u<m \}} \big(u(x)-m(x)\big)^2\, d x\leq
L\,V_m (u).
\end{align*}
The proof of inequality \eqref{sec} is in fact analogous. Due to \eqref{low-p-lap-ob-cndtn} and arguing as above  there is $L>0$ and $\delta>0$ such that, for all $x\in [0,l]$ and $s\in (0,\delta]$
\begin{equation}
s\Big(\Delta_p M(x)+f\big(x,M(x)+s\big)\Big)\leq  Ls^2.  \label{sup-sol-lip}
\end{equation}
Since $u\in B(K_M,\delta)\setminus K_M$ we get $0\leq (u-M)_+<\delta$ a.e. on $[0,l]$ Now using \eqref{sup-sol-lip} and reasoning as in (i) we get that
\begin{align*}
&  \int_{0}^{l} (u-M)_+\big(\Delta_p u + f(x,u)\big) \, dx \leq\\ &\leq -\int_{\{u>M \}} (u'-M')\big(|u'|^{p-2} u' - |M'|^{p-2}M'\big) \, dx + L\int_{\{ u>M\}}(u-M)^2 \, dx\leq L\,V_M(u),
\end{align*}
which completes the proof.
\end{proof}
As a conclusion from  Lemma \ref{p-sol-der}, Proposition \ref{p-lap-dir-der}
and  Theorem \ref{thm-inv1-prime} we get the following result.
\begin{thm}
If $m$ and $M$ satisfy \eqref{low-p-lap-ob-cndtn}, then any integral solution of \eqref{p-laplace-1-dim} such that $u(0,\cdot)=u_0\in K_m\cap K_M$ stays there for all  times $0\leq t\leq T$.\hfill $\square$
\end{thm}

\subsection{Obstacle Problem for Reaction Diffusion Equation}
Let $\Omega\subset \R^N$ be open bounded with Lipschitz boundary $\partial \Omega$ and
consider the following parabolic problem
\begin{equation}\label{reaction-diffusion}
\left\{
\begin{array}{ll}
u_t = \Delta u + f(x,u),&  x\in\Omega, \, t>0,\\
u(x,t) = 0, & x\in\partial\Omega, \, t>0,\\
\end{array}
\right.
\end{equation}
where $f:\overline{\Omega}\times \mathbb{R}\to \mathbb{R}$ is continuous and  such that $f(x,0)=0$ for all $x\in \overline{\Omega}$.
We shall deal with the  related obstacle problem, i.e. we look for solutions $u$ of  \eqref{reaction-diffusion} such that
\begin{equation}\label{obstacle-problem}
u(x,t) \geq m(x)\;\; \text{and}\;\; u(x,t) \leq M(x)\;\; \text{for all}\;\; x\in \Omega,\; t>0,
\end{equation}
where obstacles $m, M \in C^2(\Omega)$ such that
\begin{equation}\label{dot m_M}m|_{\part\Omega}\leq 0\;\; \text{and}\;\; M|_{\part\Omega}\geq 0\end{equation}
are given. As above we look for conditions on $f$, $m$ and $M$ implying that for any $u_0$ such that $m\leq u_0\leq M$ on $\Omega$ all solutions of \eqref{reaction-diffusion} starting at $u_0$ satisfy \eqref{obstacle-problem}.\\
\indent Let $X=C_0(\Omega)=\big\{ u \in C(\overline{\Omega})\mid u(x)=0\;\; \text{for all}\; x\in\partial \Omega\big\}$. Define $A:D(A)\to X$ by $Au:=-\Delta u$ (here $\Delta$ stands for the $L^2$-realization of the Laplacian), where $u\in D(A)$ and
$$
D(A):=\big\{u\in X\cap H_{0}^{1}(\Omega) \mid \Delta u \in X\big\}.
$$
It is known (see e.g. \cite[Prop. 2.6.7]{Caz-Har}) that $A$ is $m$-accretive and its domain is dense in $X$. Define $F:X\to X$ by $F(u)(x):= f\big(x,u(x)\big)$ for $u\in X$ and $x\in\overline{\Omega}$. It is clear that $F$ is well-defined and continuous. By a solution to \eqref{reaction-diffusion} on $[0,T]$ we mean a function
\begin{equation}\label{solution-regularity}
u\in C\big([0,T],X\big) \cap C\big((0,T],H_{0}^{1}(\Omega)\big) \cap C^{1}\big((0,T],L^2(\Omega)\big)
\end{equation}
such that $\Delta u(t)\in L^2(\Omega)$ and
\begin{equation}\label{abstract-reactiion-fiffusion}
\dot u (t) = \Delta u(t) +F\big(u(t)\big),\;\; \text{for each}\;\; t\in (0,T],
\end{equation}
i.e. $u$ is a strong solution to
\begin{equation}\label{A_L^2-obstacle-eq}
\dot u(t)=-A_{L^2} u(t)+F\big(u(t)\big),\;\; \text{for each}\;\; t\in (0,T],
\end{equation}
where $A_{L^2}:D(A_{L^2}) \to L^2(\Omega)$ is given by
$$
A_{L^2} u:=-\Delta u, \ D(A_{L^2}):=\{ u\in H_{0}^{1}(\Omega)\mid \Delta u\in L^2(\Omega)\}.
$$
It is well-known that $A_{L^2}$ generates a strongly continuous semigroup of contractions, i.e. it is an $m$-accretive operator with dense domain
and that the part of $A_{L^2}$ in $X$ is equal to $A$ (see e.g. \cite{Caz-Har}). Observe that, since any solution to \eqref{reaction-diffusion} in the above sense is, in particular, a strong solution to \eqref{A_L^2-obstacle-eq}, each solution of \eqref{reaction-diffusion} with $u(0)=u_0$ is also an integral solution of
\begin{equation}\label{L^2-eq_ic}
\dot u(t) = -A_{L^2}u(t)+w(t) \ \mbox{ on } [0,T], \ u(0)=u_0
\end{equation}
with $w=F \circ u$. On the other hand, if $\tilde u$ is an integral solution of
\begin{equation}\label{reaction-diff-X}
\dot u(t)=-Au(t)+F\big(u(t)\big), \ t\in [0,T],
\end{equation}
with $\tilde u(0)=u_0$, i.e. $\tilde u=u_A (\cdot,u_0,w)$, then $\tilde u=u_{A_{L^2}}(\cdot;u_0, w)$, i.e. $\tilde u$ is also an integral solution of
\eqref{L^2-eq_ic} (see the proof of Proposition \ref{reg-trick}). This means that $u=\tilde u$. Therefore, any solution of
\eqref{reaction-diffusion} is an integral solution of \eqref{reaction-diff-X}. Hence,
we are able to use Theorem \ref{thm-inv1-prime} for \eqref{reaction-diff-X} and Remark \ref{how-to-compute} (1) to verify condition \eqref{dense-prime}.\\
\indent The obstacle conditions \eqref{obstacle-problem} can be rewritten as
$$u(t)\in K_m\cap K_M,
$$
where
$$K_m :=\big\{u\in X\mid u(x)\geq m(x)\;\; \text{for}\;\; x\in \Omega\big\},\;\;\;
K_M:=\big\{u\in X\mid u(x)\leq M(x)\;\; \text{for}\;\; x\in \Omega\big\}.
$$
If we define $V_m, V_M:X\to\R$ by
$$
V_m (u):= \frac{1}{2}\int_{\Omega} (u-m)_{-}^{2} \, dx,\;\;
V_M (u):= \frac{1}{2}\int_{\Omega} (u-M)_{+}^{2} \, dx,
$$
then clearly  $K_m=\big\{u\in X\mid V_m(u)\leq 0\big\}$ and $K_M =\big\{u\in X\mid V_M(u)\leq 0\big\}$.
\begin{rem}{\em  It is well known that any solution of \eqref{reaction-diffusion} with the regularity given by \eqref{solution-regularity} is also a mild solution of \eqref{reaction-diffusion}, i.e. an integral solution to the problem
$\dot u=-Au+F(u)$ (see Remark \ref{dod info} (c)). Conversely, if  $F$ is locally Lipschitz, then a mild solution of \eqref{reaction-diffusion} satisfies \eqref{solution-regularity} (see e.g. [10, Th. 5.2.1]). Hence, for a locally Lipschitz $f$, in view of Proposition \ref{prop lip}, we infer that, for any $u_0\in X$, the problem \eqref{reaction-diffusion} admits a unique solution $u:[0,\tau_{u_0})\to X$ with the initial condition $u(0) = u_0$ and satisfying \eqref{solution-regularity} for any $0<T < \tau_{u_0}$. In the general case, when $f$ is just continuous, one has only the local existence of mild solutions (see e.g. \cite{Vrabie}).}
\end{rem}
We shall need the differentiability of $V_m$ and $V_M$ along trajectories of \eqref{reaction-diffusion}. Exactly as in Lemma \ref{p-sol-der} we get that
\begin{lm}\label{diff-const-props}
If $u:[0,T]\to X$ satisfies \eqref{solution-regularity}, then $V_m\circ u$ and $V_M \circ u$
are absolutely continuous on compact subsets of $(0,T]$ and, for almost all $t\in (0,T]$,
\begin{align*}
(V_m \circ u)'(t) =& - \int_{\Omega}(u-m)_-\Big(\Delta u(t)+f\big(x,u(t)\big)\Big) \,dx,\\
(V_M \circ u)'(t) =&  \int_{\Omega} (u-M)_+\Big(\Delta u(t)+f\big(x,u(t)\big)\Big) \,dx.
\end{align*}
\end{lm}

\begin{prop}\label{directional-der}
Suppose that
\begin{align}\label{lower-obstacle-cndtn}
\limsup_{s\to 0^-}\, \frac{\Delta m(x) + f\big(x,m(x)+s\big)}{s} < +\infty,\\
\label{lower-obstacle-cndtn1}\limsup_{s\to 0^+}\,\frac{\Delta M(x) + f\big(x,M(x)+s\big)}{s} < + \infty
\end{align}
uniformly with respect to $x\in \overline\Omega$. Then there are $C>0$ and $\delta>0$ such that, for any $u\in D(A_{L^2})\cap \big(B(K_m, \delta)\setminus K_m\big)$,
$$
-\int_{\Omega} (u-m)_-\big(\Delta u + f(x,u)\big) \, dx \leq C\,V_m (u),
$$
and, for any $u\in D(A_{L^2}) \cap\big(B(K_M,\delta)\setminus K_M\big)$,
$$
\int_{\Omega}(u-M)_+\big(\Delta u+f(x,u)\big)\, dx \leq C\,V_M (u).
$$
\end{prop}
\begin{proof}
In view of \eqref{lower-obstacle-cndtn} there is $C>0$ and $\delta>0$ such that, for any  $x\in \Omega$ and $s\in [-\delta,0)$,
$$
\frac{\Delta m(x)+f\big(x,m(x)+s\big)}{s} < C,\; \text{i.e.}\;\; s\Big(\Delta m(x) + f\big(x,m(x)+s\big)\Big) < C|s|^2.
$$
Note that in view of \eqref{dot m_M}, $(u-m)\in H_0^1(\Omega)$;  this together with Poincar\'{e}'s inequality and again \cite[Section 7.4]{Gil} gives
\begin{align*}
&-\int_{\Omega}(u-m)_-\big(\Delta u + f(x,u)\big) \, dx = \int_{\Omega}(u-m)_-\big(\Delta (u-m)  + \Delta m+ f(x,u)\big) \, dx\\
&\leq\int_{\Omega}\nabla\big[(u-m)_-\big]\nabla (u-m)\, dx -\int_{\Omega}(u-m)_-\big(\Delta m+f(x,u)\big) \, dx\\
&\leq - \int_{\Omega}\big|\nabla (u-m)_-\big|^2\,dx+\int_{\{u<m\}}(u-m)\big(\Delta m+ f(x,u)\big) \, dx\\
&\leq - 2\lambda_1 V_m (u) + C\int_{\{u<m\}} \big|u(x)-m(x)\big|^2\, d x\leq  -2\lambda_1 V_m (u) + C\,V_m (u) = (C-2\lambda_1)V_m (u),
\end{align*}
where $\lambda_1$ stands for the first eigenvalue of the (negative) Laplacian.\\
\indent Inequality \eqref{lower-obstacle-cndtn1} implies the existence of $C>0$ and $\delta>0$ such that, for all $x\in \Omega$ and $s\in (0,\delta]$
$$
\frac{\Delta M(x)+f(x,M(x)+s)}{s} < C,\;\; \text{i.e.}\;\; s\Big(\Delta M(x) + f\big(x,M(x)+s\big)\Big) < C |s|^2.
$$
Then, arguing as above, we show that
$$\int_{\Omega}(u-M)_+\big(\Delta u + f(x,u)\big) \, dx \leq (C-2\lambda_1)\,V_M(u).\eqno\qedhere$$
\end{proof}
\begin{rem} {\em  Observe that if $m$ (resp. $M$) is a subsolution (resp. supersolution) of the stationary problem
related to \eqref{reaction-diffusion}, i.e.
$$-\Delta m(x)\leq f\big(x,m(x)\big) \;\; \left(\text{resp.}\;-\Delta M(x)\geq f\big(x,M(x)\big)\right)\;\; \text{for}\;\; x\in \Omega,$$
and
$$
\limsup_{s\to 0^-} \, \frac{f\big(x,m(x)+s\big) -f\big(x,m(x)\big)}{s} < +\infty\;\; \left(\text{resp.}\;\;\limsup_{s\to 0^+} \, \frac{f\big(x,M(x)+s\big) -f\big(x,M(x)\big)}{s} < +\infty\right),
$$
uniformly with respect to $x\in \overline \Omega$, then \eqref{lower-obstacle-cndtn} (resp.  \eqref{lower-obstacle-cndtn1} holds. It is clear that each of the latter conditions is always satisfied whenever $f$ is locally Lipschitz.}
\end{rem}
Combining Lemma \ref{diff-const-props} with Proposition \ref{directional-der} yields the following.
\begin{cor}
Assume that \eqref{lower-obstacle-cndtn} and \eqref{lower-obstacle-cndtn1} hold and let $\delta>0$ and $C>0$ be as in Proposition \ref{directional-der}. If $u:[0,T]\to X$ is a solution to \eqref{reaction-diffusion}, then
$$
(V_m \circ u)'(t) \leq C\,V_m \big(u(t)\big), \ \ \mbox{ for a.e. } \ t\in (0,T],
\ \mbox{  whenever } \ u(t)\in B(K_m,\delta) \setminus K_m,
$$
and
$$
(V_M \circ u)'(t) \leq C\,V_M\big(u(t)\big), \ \ \mbox{ for a.e.} \ t\in (0,T],
\ \mbox{ whenever } \ u(t)\in B(K_M,\delta) \setminus K_M. \eqno\square
$$
\end{cor}
This together with Theorem \ref{thm-inv1-prime} allows us to conclude that, under our assumptions on $f$, $m$ and/or $M$, if a solution of the reaction diffusion problem \eqref{reaction-diffusion} starts above the obstacle $m$ and below the obstacle $M$ then it does not cross the obstacle(s) in positive times. Formally we state it as the following.
\begin{thm}
If $f$ satisfies conditions \eqref{lower-obstacle-cndtn} and \eqref{lower-obstacle-cndtn1}, and $u$ is a solution of \eqref{reaction-diffusion} such that \eqref{obstacle-problem} holds for $t=0$, then \eqref{obstacle-problem} holds for all positive times $t$.\hfill $\square$
\end{thm}

\subsection{Age-Structured Population Model}
Consider the following first order partial differential problem
\begin{equation}\label{age-dep-model}
\left\{
\begin{array} {l}
u_t(x,t) + u_x(x,t) = f\big(x,u(x,t)\big),\;\; x\in (0,a),\;t>0,\\
u(0,t) = \displaystyle{\int_{0}^{a}}\beta (x)u(x,t)\,dx,\;\; t>0.
\end{array}
\right.
\end{equation}
This is McKendrick's model for an age dependent population with the so-called birth function $\beta:[0,a]\to [0,+\infty)$ satisfying the condition
\begin{equation}\label{birth-cond}
\int_{0}^{a} \beta(x)\,dx <1,
\end{equation}
and the so-called "migration" component $f: [0,a ] \times \mathbb{R}\to \mathbb{R}$ that is locally Lipschitz in the second variable uniformly with respect to the first one, i.e. for any $s\in \mathbb{R}$ there exist $L>0$ and $\delta>0$  such that, for any $x\in [0,a]$ and $s_1, s_2\in [s-\delta, s+\delta]$ one has
$\big|f(x,s_1)-f(x,s_2)\big|\leq L|s_1-s_2|$.\\
\indent We provide additional conditions on $f$ and $m:[0,a]\to\mathbb{R}$ assuring that any solution $u$ of \eqref{age-dep-model} such that
\begin{equation} \label{m-obst-age-0}
u(x,0) \geq m(x)\;\; \text{for all}\;\; x\in [0,a]
\end{equation}
will preserve the inequality for positive times $t$, i.e.
\begin{equation}\label{m-obst-age-t}
u(x,t)\geq m(x)\;\; \text{for all}\;\; x\in [0,a] \text{ and } t\geq 0 .
\end{equation}
\indent Let $X=C[0,a]$ and  define  $A:D(A) \to X$ by
$$
Au:= u',\;\; u\in D(A):=
\left\{u\in X\cap C^1(0,a)\mid u(0)=\int_{0}^{a} \beta(x) u(x)\, dx \right\}.
$$
Condition \eqref{birth-cond} implies, by a straightforward calculation, that $A$ is $m$-accretive. Let $F:X\to X$ be defined by $\big[F(u)\big](x):=f\big(x,u(x)\big)$ for $u\in X$ and $x\in [0,a]$. Problem \eqref{age-dep-model} can be rewritten as
\begin{equation}\label{abst-age-dep}
\dot u(t) = -Au(t)+F\big(u(t)\big),\;\; t\in [0,\tau).
\end{equation}
Define also $A_{L^2}:D(A_{L^2}) \to L^2[0,a]$ by
$$A_{L^2}u:=u'\;\; \text{for}\;\; u\in D(A_{L^2}):=\left\{
u\in H^1 (0,a)\mid u(0) = \int_{0}^{a}\beta(x)u(x)\,dx \right\}.
$$
For $\lambda:= \frac{1}{2}\|\beta\|_{L^2}^2$, by the use of the Cauchy-Schwarz inequality, one has
\begin{align*}
\la u, A_{L^2} u+\lma u\ra_{L^2} & = \int_{0}^{a} u u' \, dx+\lambda\|u\|^2_{L^2}  =\frac{1}{2}u^2(a)-\frac{1}{2}\left(\int_0^a\beta u\,dx\right)^2+\lambda\|u\|_{L^2}^2\\
&\geq \frac{1}{2}u^2(a)-\frac{1}{2}\|\beta\|_{L^2}^2\|u\|^2_{L^2}+\lambda\|u\|_{L^2}^2
=\frac{1}{2}u^2(a)\geq 0.
\end{align*}
Hence $A_{L^2}+\lambda I$ is accretive. The $m$-accretivity of $A_{L^2}+\lambda I$ may be proved in the same  way as for $A$. Therefore $A_{L^2}$ is $\lambda$-$m$-accretive.\\
\indent In view of Proposition \ref{reg-trick} we get the following regularity result.
\begin{lm} \label{ad-reg-sol}
For any integral solution $u:[0,\tau)\to X$ of \eqref{abst-age-dep} with $u(0)\in D(A)$  one has $u\in W^{1,1}([0,T],L^2[0,a])$, for any $T\in (0,\tau)$, $u(t)\in D(A_{L^2})$ and
\begin{equation*} 
\dot u (t) = - A_{L^2}u(t) + F\big(u(t)\big),\;\; \text{for a.e.}\;\; t\in [0,\tau).
\end{equation*}
\end{lm}
\noindent This means that we are able to apply the idea from Remark \ref{how-to-compute} (2) in order to verify the assumption (iii)' from Theorem \ref{thm-inv2-prime}.\\
\indent Define $V:X\to\R$ by
$$
V(u):=\frac{1}{2} \int_{0}^{a}(u-m)_-^2 \, dx,\;\; u\in X.
$$
Let
$$
K:=\big\{u\in X\mid u(x)\geq m(x)\;\; \text{for}\;\; x\in [0,a]\big\}.
$$
\begin{prop}\label{age-dep-der}
If $m\in C^1(0,a)\cap C[0,a]$,
\begin{equation}\label{ad-tang}
- m'(x)+f\big(x,m(x)\big) \geq 0\;\; \text{for all}\;\; x\in [0, a],
\end{equation}
and
\begin{equation}\label{ad-tang1}
m(0)=\int_{0}^{a} \beta(x) m(x) \, dx,
\end{equation}
then there exist $C>0$ and $\delta>0$ such that, for any mild solution $u:[0,\tau)\to X$ of \eqref{abst-age-dep}, the function
$V\circ u$ is a.e. differentiable and, for a.e. $t\in [0,\tau)$,
$$
(V\circ u)'(t) \leq C\,V\big(u(t)\big)\;\; \text{if}\;\;  u(t)\in B(K,\delta)\setminus K.
$$
\end{prop}
\begin{proof} By \eqref{ad-tang} and the local Lipschitz condition for $f$, there are $L>0$ and $\delta>0$ such that for all $s\in (-\delta,0)$ and $x\in [0,a]$,
\begin{equation}\label{de}\begin{split}
s\Big(-m'(x)+f\big(x,m(x)+s\big)\Big)= &\; s\Big(-m'(x)+f(x,m(x)\big)\big)\\
&+ s\Big(f\big(x,m(x)+s\big)-f\big(x,m(x)\big)\Big)\leq L|s|^2.
\end{split}
\end{equation}
According to Lemma \ref{ad-reg-sol}, $u$ is almost everywhere differentiable as a function into $L^2[0,a]$. If we define $\wt V:L^2[0,a]\to \R$ by
$$
\wt V(u):=\frac{1}{2}\int_{0}^{a}(u-m)_-^2\,dx,\;\; u\in L^2[0,a],
$$
then, by the chain rule,  $u(t)\in D(A_{L^2})$ for a.e. $t\in [0,\tau)$,  and
$$
(V\circ u)'(t)  =  (\wt V\circ u)'(t) = \Big[\wt V'\big(u(t)\big)\Big]\big(\dot u(t)\big) =
-\int_{0}^{a} u(t)_-\big(-A_{L^2}u(t) + F\big(u(t)\big)\big)\, dx.
$$
Put $\varphi:=u(t)$ and observe that $-\delta<\vp(x)-m(x)\leq 0$ if $x\in\{y\in [0,a]\mid (\vp-m)_-(y)\neq 0\}$. Hence
\begin{align*}
-(V\circ u)'(t) & = \int_{0}^{a}(\varphi-m)_-\big(-\varphi' + f(x,\varphi)\big)\,dx\\
& =\int_{0}^{a}(\varphi-m)_-\big[(\varphi-m)'\big]\, dx+\int_0^a(\vp-m)_-\big(-m'+f(x,\vp)\big)\,dx \\
&= \int_{0}^{a}(\varphi-m)_-\big[(\vp-m)_-\big]'\,dx-\int_0^a(\vp-m)\big(-m'+f(x,\vp)\big)\,dx\\
&\geq \frac{1}{2}\left(\big[(\vp-m)_-\big]^2(a)-\big[(\vp-m)_-\big]^2(0)\right)-
L\int_{\{\vp<m\}}(\vp-m)^2\,dx\\
&\geq -\frac{1}{2}\big[(\vp-m)_-\big]^2(0)-L\int_{\{\vp<m\}}(\vp-m)^2 \, dx.\end{align*}
Note that, since $\varphi,m\in D(A_{L^2})$ and $\beta\geq 0$, one gets
$$
(\varphi-m)(0) = \int_{0}^{a}\beta(\varphi-m)\, dx = \int_{0}^{a}\beta(\varphi-m)_+\, dx -\int_{0}^{a}\beta (\varphi-m)_-\, dx,
$$
i.e.
$$
-(\varphi-m)(0) \leq \int_{0}^{a} \beta(\varphi-m)_-\, dx\;\; \text{and thus}\;\; (\vp-m)_-(0)\leq\int_0^a\beta(\vp-m)_-\,dx.
$$
By the Cauchy-Schwarz inequality we get
$$
(V\circ u)'(t)  \leq   \frac{1}{2}\left(\int_{0}^{a} \beta(\varphi-m)_-\, dx \right)^2  + LV(\varphi) \leq \left(\frac{1}{2}\|\beta\|_{L^2}^2+L\right)V(\varphi).\eqno\qedhere
$$
\end{proof}
Combining Proposition \ref{age-dep-der} with Theorem \ref{thm-inv2-prime} we conclude with the following.
\begin{thm}
If $f$ satisfies \eqref{ad-tang} and \eqref{ad-tang1}, then any mild solution $u$ of \eqref{age-dep-model} satisfying condition \eqref{m-obst-age-0} satisfies \eqref{m-obst-age-t} for all positive times $t$.
\end{thm}

\end{document}